\documentclass[a4paper,final]{amsart}
\usepackage{amsfonts}
\usepackage{amssymb}
\usepackage{amsfonts}
\usepackage{amssymb}
\usepackage{amsthm}
\usepackage{amsmath}
\usepackage{a4wide}

\theoremstyle{plain}

\newtheorem{teo}{Theorem}[section]
\newtheorem{lemma}[teo]{Lemma}
\newtheorem{prop}[teo]{Proposition}
\newtheorem{cor}[teo]{Corollary}

\theoremstyle{definition}

\newtheorem{akn}[]{}

\theoremstyle{remark}

\newtheorem{rem}[teo]{Remark}

\numberwithin{equation}{section}

\def\div{\operatornamewithlimits{div}\nolimits}

\def\R{{{\mathbb R}}}
\def\SS{{{\mathbb S}}}

\def\NN{{{\mathbb N}}}

\def\composed{\circ}
\def\comp{\composed}

\def\tr{\text{\rm {Trace} } }
\def\HHH{{\mathrm H}}

\def\Ri{\R^{n+1}}
\def\AAA{{\mathrm A}}
\def\bbigstar{\operatornamewithlimits{\text{\Large{$\circledast$}}}}
\def\EEE{{\mathrm E}}
\def\pol{{\mathfrak{p}}}
\def\qol{{\mathfrak{q}}}
\def\RRR{{\mathrm R}}
\def\Ric{{\mathrm {Ric}}}
\def\Scal{{R}}

\begin{document}

\title[Smooth Geometric Evolutions]{Smooth Geometric Evolutions of
Hypersurfaces}

                                \author[Carlo Mantegazza]{Carlo Mantegazza}
                                \address[Carlo Mantegazza]{Scuola
                                Normale Superiore, Pisa, 56126, Italy}
                                \email[C. Mantegazza]{mantegaz@sns.it}

\date{\today}

\begin{abstract} We consider the gradient flow associated to the
following functionals
$$
{\mathcal F}_m(\varphi)=\int_M1+\vert\nabla^m\nu\vert^2\,d\mu\,.
$$
The functionals are defined on hypersurfaces immersed in $\Ri$ via
a map $\varphi:M\to\Ri$, where $M$ is a smooth closed and 
connected $n$--dimensional manifold without boundary.\\
Here $\mu$ and $\nabla$ are respectively the canonical measure and
the Levi--Civita connection on the Riemannian manifold $(M,g)$, where
the metric $g$ is obtained by pulling back on $M$ the usual metric of
$\Ri$ with the map $\varphi$. The symbol $\nabla^m$ denotes the
$m$--th iterated covariant derivative and $\nu$ is a unit normal local
vector field to the hypersurface.

Our main result is that if the order of derivation $m\in\NN$ is
strictly larger than the integer part of $n/2$ then singularities in
finite time cannot occur during the evolution.

These geometric functionals are related to similar ones proposed by 
Ennio~De~Giorgi, who conjectured for them an analogous regularity 
result. In the final section we discuss the original conjecture 
of De~Giorgi and some related problems.
\end{abstract}

\maketitle


\section{Introduction}

In one of his last papers Ennio~De~Giorgi 
conjectured that any compact $n$--dimensional hypersurface in $\Ri$, 
evolving by the gradient flow of certain functionals depending on
sufficiently high derivatives of the curvature does not develop
singularities during the flow (\cite{degio5,degio6}, Sec.~5, Conj.~2 
--- see~\cite{degio6} for an English translation).\\
This result is central in his program to approximate singular
geometric flows with sequences of smooth ones.

Representing hypersurfaces in $\Ri$ as immersions $\varphi:M\to\Ri$,
we consider the gradient flow associated to the following
functionals
$$
{\mathcal F}_m(\varphi)=\int_M1+\vert\nabla^m\nu\vert^2\,d\mu
$$
where $\mu$ and $\nabla$ are respectively the canonical measure and
the Levi--Civita connection on the Riemannian manifold $(M,g)$, where
the metric $g$ is obtained by pulling back on $M$ the usual metric of
$\Ri$ via $\varphi$. We denote with $\nabla^m$ the $m$--th
iterated covariant derivative and with $\nu$ a unit normal local
vector field to the hypersurface. Finally, $\AAA$ and $\HHH$ are
respectively the second fundamental form and the mean curvature of
the hypersurface.\\
These functionals are strictly related to the ones proposed by 
De~Giorgi since, roughly speaking, the derivative of the normal is the
curvature of $M$. Though not exactly the same, they can play the
same role in the approximation process he suggested. In the end of the
paper we discuss some other possible functionals and, in particular,
the original De~Giorgi conjecture.

Our main result is that if the order of derivation $m\in\NN$ is
strictly larger than $\left[\frac{n}{2}\right]$ (where
$\left[\frac{n}{2}\right]$ denotes the integer part of $n/2$), then
singularities cannot occur.

The simplest case $n=1$ and $m=1$ is concerned with curves in the
plane evolving by the gradient flow of
\begin{equation}\label{pden1}
{\mathcal F}_1(\gamma)=\int_{\SS^1} 1+k^2\,ds
\end{equation}
since the curvature $k$ of a curve $\gamma:\SS^1\to\R^2$ satisfies
$k^2=\vert\nabla\nu\vert^2$. The global regularity in this case was
showed by Polden in the papers~\cite{polden1,polden2} which have been
a starting point for our work. 
Wen in~\cite{wen2} found results similar to Polden's ones, in
considering the flow for $\int_{\SS^1} k^2\,ds$ of curves with a
fixed length.

The very first step in attacking our problem is an analysis of 
the first variation of the functionals ${\mathcal F}_m$, which gives 
rise to a quasilinear system of partial differential equations on the
manifold  $M$.\\
The small time existence and uniqueness of a smooth flow is a
particular case of a very general result of Polden proven
in~\cite{huiskpold,polden2}. Then the long time existence is
guaranteed as soon we have suitable a priori estimates on the flow.

In the study of the {\em mean curvature flow} of a hypersurface
$\varphi:M\times[0,T)\to\Ri$, 
$$
\frac{\partial\varphi}{\partial t} = - \HHH \nu = \Delta_t\varphi\,,
$$
(which is of second order) via techniques such as varifolds, level
sets, viscosity solutions
(see~\cite{altawa,ambson,brakke,es,ilman1}),  
the maximum principle is the key tool to get comparison results and
estimates on solutions. 
In our case, even if $m=1$, the first variation and hence the corresponding
parabolic problem turns out to be of order higher than two, precisely 
of order $2m+2$, so we have to deal with equations of fourth order
at least. This fact has the relevant consequence that we cannot employ the
maximum principle to get pointwise estimates and to compare two
solutions, thus losing a whole bunch of geometric results holding 
for the mean curvature flow. In particular, we cannot expect that
an initially embedded hypersurface remains
embedded during the flow, since self--intersections can appear in
finite time (an example is given by Giga and Ito in~\cite{gigaito1}).
By these reasons, techniques based on the description of the
hypersurfaces as level sets of functions seems of difficult
application in this case and therefore we adopt a parametric
approach as in the work of Huisken~\cite{huisk1}.

Despite the large literature on the mean curvature flow, fourth or
even higher order flows appeared only recently. Besides the cited
works of Polden and Wen, we quote the work of
Escher, Mayer and Simonett~\cite{escmaysim} on the {\em surface
  diffusion flow} (see also the references therein) 
$$
\frac{\partial\varphi}{\partial t}=( \Delta_t\HHH )\nu
$$ 
and of Simonett~\cite{simonett1} on the gradient flow of the {\em
  Willmore functional} (see~\cite{willmore}) 
$$
{\mathcal W}(\varphi)=\int_M \vert\AAA\vert^2\,d\mu  
$$
defined on surfaces immersed in $\R^3$. In these papers it is shown
the long term existence and convergence of the flow for initial data
which are $C^{2,\alpha}$--close to a sphere.\\
In the article of Chru\'sciel~\cite{chrusciel2}, the global existence of
a fourth order flow of metrics on a two--dimensional Riemannian manifold is
applied to construct solutions of Einstein vacuum equations
representing an isolated gravitational system, called
Robinson--Trautman metrics.\\
Another problem considered by Polden in~\cite{polden2,polden3} 
is the conformal evolution of a metric $g$ on a two--dimensional
manifold $M$ by the gradient flow of the functional
$$
{\mathcal R}(g)=\int_M F(\Scal)\,d\mu
$$
where $\Scal$ is the scalar curvature of $(M,g)$ and $F$ is an even,
smooth and strictly convex function.\\
Finally, in a very recent paper~\cite{kuschat1} Kuwert and Sch\"atzle 
study the global existence and regularity of the gradient flow of the 
{W}illmore functional for general initial data.\\
Our work borrows from~\cite{chrusciel2,polden1,polden2,polden3} the
basic idea of using interpolation inequalities as a tool to get a
priori estimates.

We want to remark here that a strong motivation for the study of
these flows is the fact that, in general, regularity is not shared by 
second order flows, with the notable exceptions of the evolution by
mean curvature of embedded curves in the plane
(see~\cite{gaha1,gray1,huisk2}) and of convex
hypersurfaces (see~\cite{huisk1}). So our result opens the possibility
to approximate  canonically singular flows with smooth ones by singular
perturbation arguments (see~\cite{degio5,degio6} and
Section~\ref{remmmm}).

In order to show regularity, a good substitute of the pointwise
estimates coming from the maximum  principle, are suitable 
estimates on the second fundamental form in Sobolev spaces,
using Gagliardo--Nirenberg interpolation type inequalities for
tensors. Since the constants involved in these inequalities 
depends on the Sobolev constants and these latter on the
geometry of the hypersurface where the tensors are defined, before
doing estimates we absolutely need some uniform control independent of 
time on these constants. In~\cite{polden1} these controls are obvious as the
constants depend only on the length, on the contrary, much
more work is needed in~\cite{chrusciel2,kuschat1,polden3}, because of
the richer geometry of surfaces.

In our case, we will see that if $m$ is large enough, the functional
${\mathcal F}_m$, which decreases during  the flow, controls the
$L^p$ norm of the second fundamental form for some exponent $p$ larger
than the dimension. This fact, combined with a universal Sobolev type
inequality due to Michael  and Simon~\cite{micsim}, where the dependence
of the constants on the curvature is made explicit, 
allows us to get an uniform bound on the Sobolev 
constants of the evolving hypersurfaces and then to obtain 
time--independent estimates on curvature and all its derivatives in
$L^2$. These bounds imply in turn the desired pointwise estimates and 
the long time existence and regularity of the flow.

In the last section we will discuss some possible extensions of our
results, some open problems and the related conjectures of De~Giorgi.

\smallskip

\begin{akn} We are grateful to Gerhard Huisken for many discussions 
  about geometric flows during his visit at the Scuola Normale
  Superiore of Pisa. Moreover, we wish to thank Luigi Ambrosio for 
his constant  encouragement and invaluable help in several occasions.

Our work would have been impossible without the enlightening
  mathematical insight of Ennio~De~Giorgi. This paper is dedicated to
  his memory.
\end{akn}

\section{Notation and Preliminaries}\label{prelim}

We devote this section to introduce the basic notations and 
facts about  differentiable and Riemannian manifolds we need in the
paper, a good reference for this introduction is~\cite{gahula} or the
first part of~\cite{petersen1}.

The main objects of the paper are $n$--dimensional closed hypersurfaces
immersed in $\Ri$, that is, pairs $(M,\varphi)$ where $M$ is an
$n$--dimensional smooth manifold, compact, connected with empty
boundary, and a smooth map  $\varphi:M\to\R^{n+1}$ such that the rank
of $d\varphi$ is everywhere equal to $n$.

The manifold $M$ gets in a natural way a metric tensor $g$ turning it 
in a Riemannian manifold $(M,g)$, by pulling back the standard scalar
product of $\Ri$ with the immersion map $\varphi$.

Taking local coordinates around $p\in M$ given by a chart
$F:\R^n\supset U\to M$,  we identify the map $\varphi$ with its
expression in coordinates $\varphi\comp F:\R^n\supset U\to\Ri$, then
we have local basis of $T_{p}M$ and $T^*_{p}M$, respectively given 
by vectors $\left\{\frac{\partial}{\partial x_i}\right\}$ and
covectors $\{dx_j\}$.

We will denote vectors on $M$ by $X=X^i$, which means 
$X=X^i\frac{\partial}{\partial x_i}$, covectors by
$Y=Y_j$, that is, $Y=Y_jdx_j$ and a general mixed tensor with
$T=T^{i_1\dots i_k}_{j_1\dots j_l}$, where the indices refer to the
local basis.

Sometimes we will need also to consider tensors along
$M$, viewing it as a submanifold of $\Ri$ via the map $\varphi$, in
that case we will use the Greek indices to denote the components of
such tensors in the canonical basis $\{e_\alpha\}$ of $\Ri$, for
instance, given a vector field $X$ along $M$, not necessarily tangent,
we will have  $X=X^\alpha e_\alpha$.

{\em In all the paper the convention to sum over repeated 
indices will be adopted.}

The inner product on $M$, extended to tensors, is given by
$$
g(T , S)=g_{i_1s_1}\dots g_{i_k s_k}g^{j_1
z_1}\dots g^{j_l z_l}  T^{i_1\dots
  i_k}_{j_1\dots j_l}S^{s_1\dots
  s_k}_{z_1\dots z_l}
$$
where $g_{ij}$ is the matrix of the coefficients of the metric
tensor in the local coordinates and $g^{ij}$ is its 
inverse. Clearly, the norm of a tensor is
$$
\vert T\vert=\sqrt{g(T , T)}\,.
$$

The scalar product in $\Ri$ will be denoted with 
$\langle\cdot\,\vert\,\cdot\rangle$. As the metric $g$ is obtained
pulling it back with $\varphi$, we have 
$$
g_{ij}(x)=\left\langle\frac{\partial{{\varphi}}(x)}{\partial
    x_i}\,\biggl\vert\,\frac{\partial{{\varphi}}(x)}{\partial
    x_j}\right\rangle\,.
$$

The canonical measure induced by the metric $g$ is given by
$\mu=\sqrt{G}\,{\mathcal L}^n$ where $G=\det(g_{ij})$ and
${\mathcal L}^n$ is the standard Lebesgue measure on $\R^n$.

The {\em second fundamental form} $\AAA= h_{ij}$ of $M$ is the 2--tensor defined
as follows:
$$
h_{ij}(x)=-\left\langle\nu(x)\,\biggl\vert\,\frac{\partial^2{{\varphi}}(x)}
{\partial x_i\partial x_j}\right\rangle\,,
$$
the {\em mean curvature} $\HHH$ is the trace of $\AAA$, 
\begin{equation}\label{meanc}
\HHH(x)=g^{ij}(x)h_{ij}(x)\,.
\end{equation}
The induced covariant derivative on $(M,g)$ of a vector field $X$ is
given by
$$
\nabla_jX^i=\frac{\partial}{\partial
  x_j}X^i+\Gamma^{i}_{jk}X^k
$$
where the Christoffel symbols $\Gamma=\Gamma^{i}_{jk}$ are expressed
by the following formula,
$$
\Gamma^{i}_{jk}=\frac{1}{2} g^{il}\left(\frac{\partial}{\partial
    x_j}g_{kl}+\frac{\partial}{\partial
    x_k}g_{jl}-\frac{\partial}{\partial
    x_l}g_{jk}\right)\,.
$$
In all the paper the covariant derivative $\nabla T$ of a tensor
$T=T^{i_1\dots i_k}_{j_1\dots j_l}$ will be denoted by 
$\nabla_sT^{i_1\dots i_k}_{j_1\dots j_l}=(\nabla T)^{i_1\dots
  i_k}_{sj_1\dots j_l}$. 

With $\nabla^m T$ we will mean the $k$--th iterated covariant
derivative of a tensor $T$.

We recall that the gradient $\nabla f$ of a function and the 
divergence $\div X$ of a vector field at a point $p\in(M,g)$ are defined
respectively by
$$
g(\nabla f(p) , v)=df_p(v)\qquad\forall v\in T_pM
$$
and
$$
\div X=\tr \nabla X=\nabla_iX^i=\frac{\partial}{\partial
  x_i}X^i+\Gamma^{i}_{ik}X^k\,.
$$
Notice that considering $M$ as a submanifold of $\Ri$, 
if $\{e_i\}\in\Ri$ is an orthonormal basis of $T_pM$ we can 
express the divergence of $X$ as
$$
\div X(p) = g(e_i , \nabla_{e_i}X) = 
\langle e_i\,\vert\,\nabla^M_{e_i}X\rangle = 
\nabla^M_{e_i}\langle X\,\vert\,e_i\rangle = 
\nabla_{e_i}\langle X\,\vert\,e_i\rangle
$$
where $\nabla^M$ denotes the projection on $T_pM$ of the covariant
derivative of $\Ri$.\\
Using this last expression we can define the divergence of a general, not
necessarily tangent, vector field $X$ along $M$ as a 
Riemannian submanifold of $\Ri$.\\
Such definition is useful in view of the following {\em tangential
  divergence formula} (see~\cite{simon}, Chap.~2, Sec.~7), 
\begin{equation}\label{divtang}
\int_M\div X\,d\mu=\int_M \langle\nu\,\vert\, X\rangle \HHH \,d\mu
\end{equation}
holding for every vector field $X$ along $M$.\\
Notice that the right term is well defined since, by
definition~\eqref{meanc}, $\HHH \, \nu$ is independent of the choice
of the local unit normal $\nu$. Moreover, if $X$ is a tangent vector
field we recover the usual {\em divergence theorem}
$$
\int_M\div X\,d\mu=0\,.
$$

The Laplacian $\Delta T$ of a tensor $T$ is
$$
\Delta T=g^{ij}\nabla_i\nabla_jT\,.
$$

The Riemann tensor, the Ricci tensor and the scalar curvature are
expressible via the second fundamental form as follows,
\begin{align*}
\RRR_{ijkl}\,=&\,h_{ik}h_{jl}-h_{il}h_{jk}\,,\\
\Ric_{ij}\,=&\,\HHH \, h_{ij}-h_{il}g^{lk}h_{kj}\,,\\
\Scal\,=&\,\HHH^2-\vert\AAA\vert^2\,.
\end{align*}
Hence, the formulas for the interchange of covariant derivatives,
which involve the Riemann tensor, become
$$
\nabla_i\nabla_jX^s-\nabla_j\nabla_iX^s=\RRR_{ijkl}g^{ks}X^l=\RRR_{ijl}^sX^l=
\left(h_{ik}h_{jl}-h_{il}h_{jk}\right)g^{ks}X^l\,,
$$
\begin{equation}\label{ichange}
\nabla_i\nabla_jY_k-\nabla_j\nabla_iY_k=\RRR_{ijkl}g^{ls}Y_s=\RRR_{ijk}^sY_s=
\left(h_{ik}h_{jl}-h_{il}h_{jk}\right)g^{ls}Y_s\,.
\end{equation}
The Codazzi equations
$$
\nabla_ih_{jk}=\nabla_jh_{ik}=\nabla_kh_{ij}
$$ 
imply the following identity (see~\cite{simons}) 
which will be crucial in the sequel,
\begin{equation}\label{codaz}
\Delta h_{ij}=\nabla_i\nabla_j\HHH + \HHH \, 
h_{il}g^{ls}h_{sj}-\vert\AAA\vert^2h_{ij}\,.
\end{equation}
Also fundamental will be the Gauss--Weingarten relations
\begin{equation}\label{gwein}
\frac{\partial^2\varphi}{\partial x_i\partial
  x_j}=\Gamma_{ij}^k\frac{\partial\varphi}{\partial
  x_k}-h_{ij}\nu\,,\qquad\frac{\partial}{\partial x_j}\nu=
h_{jl}g^{ls}\frac{\partial\varphi}{\partial x_s}\,,
\end{equation}
which easily imply $\vert\nabla\nu\vert=\vert\AAA\vert$.

Now we introduce some non standard notation 
which will be useful for  the computations of the following 
sections.

In all the paper we will write $T * S$, following
Hamilton~\cite{hamilton1}, to denote a tensor formed by contraction on 
some indices of the tensors $T$ and $S$ using the coefficients 
$g^{ij}$.\\
Abusing a little the notation, if $T_1, \dots , T_l$ is a finite
family of tensors (here $l$ is not an index of the tensor $T$), 
with the symbol
$$
\bbigstar_{i=1}^l T_i
$$
we will mean $T_1 * T_2 * \dots * T_l$\,.

We will use the symbol $\pol_{s}(T_1, \dots , T_l)$ for a polynomial
in the tensors $T_1,\dots , T_l$ and their iterated covariant
derivatives with the $*$ product like
$$
\pol_{s}(T_1, \dots , T_l)=\sum_{i_1+\dots+i_l=s}
c_{i_1\dots i_l}\,\nabla^{i_1} T_1 * \dots * \nabla^{i_l} T_{l}\,,
$$
where the $c_{i_1\dots i_l}$ are some real constants.\\
Notice that every tensor $T_{i}$ must be present in every additive 
term of $\pol_{s}(T_1, \dots , T_l)$ and there are not repetitions.\\
We will use instead the symbol $\qol^s$ when the tensors involved are
all $\AAA$ or $\nabla\nu$, repetitions are allowed and in every
additive term of there must be present every argument of
$\qol^s$, for instance,
$$
\qol^s(\nabla\nu, \AAA)=
\sum\left(\bbigstar_{k=1}^N\nabla^{i_k}(\nabla\nu)\,\,\bbigstar_{l=1}^M\nabla^{j_l}\AAA\,\right) 
\qquad\text{ with $N$, $M\geq1$.}
$$
The order $s$ denotes the sum
$$
s=\sum_{k=1}^N (i_k+1)\,+\,\sum_{l=1}^M (j_l+1)\,.
$$

\begin{rem}
Supposing that $\qol^s$ is completely contracted, that is, there are
no free indices and we get a function, then the order $s$ has the
following strong geometric meaning: if we consider the family of
homothetic immersions $\lambda\varphi:M\to\Ri$ for $\lambda>0$, they
have associated normal $\nu^\lambda$, metric $g^\lambda$, connection 
$\nabla^\lambda$ and second form $\AAA^\lambda$ satisfying the 
following rescaling equations,
$$
(\nabla^\lambda)^i\nu^\lambda=\nabla^i\nu \qquad 
(\nabla^\lambda)^j\AAA^\lambda=\lambda\nabla^j\AAA\,,
$$
$$
(g^\lambda)_{ij}=\lambda^{2}g_{ij} \qquad
(g^\lambda)^{ij}=\lambda^{-2}g^{ij}\,.
$$
Then every completely contracted polynomial $\qol^s$ in $\nabla\nu$
and $\AAA$ will have the form 
$$
\sum(\nabla^{i_1}\nabla\nu) \dots (\nabla^{i_k}\nabla\nu) \dots
(\nabla^{i_N}\nabla\nu) \, 
\nabla^{j_1}\AAA \dots \nabla^{j_l}\AAA \dots \nabla^{j_M}\AAA \, 
g^{w_1z_1} \dots g^{w_tz_t}
$$
with 
$$
s=\sum_{k=1}^N (i_k+1)\,+\,\sum_{l=1}^M (j_l+1)
$$
and since the contraction is total it must be
$$
t=\frac{1}{2}\left(\sum_{k=1}^N (i_k+1)+\sum_{l=1}^M (j_l+2)\,\right)
$$
as the sum between the large brackets give the number of covariant
indices in the product above.\\
By this argument and the rescaling equations above, we see that
$\qol^s$ rescales as
\begin{align*}
\qol^{s}(\nabla^\lambda\nu^\lambda, \dots , \AAA^\lambda)
\,=&\,\lambda^{M-2t}\qol^{s}(\nabla\nu, \dots , \AAA)\\
\,=&\lambda^{-\left(\sum_{k=1}^N (i_k+1)+\sum_{l=1}^M
(j_l+1)\,\right)}\qol^{s}(\nabla\nu, \dots , \AAA)\\
\,=&\lambda^{-s}\qol^{s}(\nabla\nu, \dots , \AAA)\,.
\end{align*}
By this reason, with a little misuse of language, also when $\qol^s$
is not completely contracted, we will say that $s$ is the {\em
  rescaling order} of $\qol^s$.
\end{rem}

In most of the following computations only the rescaling order and the 
arguments of the polynomials involved  will be important, so we will
avoid to make explicit their inner structure.\\
An example in this spirit, are the following substitutions that we
will often apply
$$
\nabla\pol_s(T_1, \dots , T_l)=\pol_{s+1}(T_1,
\dots , T_l)\quad\text{ and }\quad\nabla\qol^z(\nabla\nu , \dots ,
\AAA)=\qol^{z+1}(\nabla\nu, \dots , \AAA)\,.
$$

We advise the reader that the polynomials $\pol_s$ and $\qol^z$ could
vary from a line to another in a computation by addition of terms
with the same rescaling order. Moreover, also the constants could vary
between different formulas and from a line to another.

\section{First Variation}\label{firstvariat}

Given an immersion $\varphi:M\to\Ri$ of a smooth closed hypersurface
in $\Ri$, we consider the following functionals for $m\geq1$,
$$
{\mathcal F}_m(\varphi)=\int_M1+\vert\nabla^m\nu\vert^2\,d\mu
$$
where $\nu$ is a local unit normal vector field to $M$ and
$\vert\nabla^m\nu\vert^2$ means
$\sum_{\alpha=1}^{n+1}\vert\nabla^m\nu^\alpha\vert^2$. The norm 
$\vert\,\cdot\,\vert$, the connection $\nabla$ and  the measure $\mu$
are all relative to the Riemannian metric $g$ which is induced on $M$
by $\Ri$ via the immersion $\varphi$. Notice that these functionals are 
well defined also without a global unit normal vector field, i.~e., 
$M$ is not orientable, because of the modulus.

In this section we are going to analyze the first variation of these
functionals. Actually, computing the exact form can be quite long but
for our purposes we need only to study some properties of its
structure.

Suppose that we have a one parameter family ${\mathcal I}$ of
immersions $\varphi_t:M\to\Ri$, with $\varphi_0=\varphi$, we compute
\begin{equation}\label{derivaz}
\delta{\mathcal F}_m(\varphi)(\mathcal
I)=\left.\frac{d}{dt}{{\mathcal F}}_m(\varphi_t)\,\right\vert_{t=0}=
\left.\frac{d}{dt}\int_M1+\vert
    \nabla^m\nu\vert^2\,d\mu_t\,\right\vert_{t=0}
\end{equation}
where clearly the metric $g$, the covariant derivative $\nabla$ and the
normal $\nu$ depend on $t$.\\
Setting $X(p)=\left.\frac{\partial}{\partial
    t}\varphi_t(p)\,\right\vert_{t=0}$ we obtain a vector field along
$M$ as a submanifold of $\Ri$ via $\varphi$. It is well known that
$$
\left.\frac{\partial}{\partial
    t}\mu_t\,\right\vert_{t=0}=\,\HHH \langle\nu\,\vert\,X\rangle\,\mu
$$
so it follows,
\begin{align*}
\left.\frac{d}{dt}{{\mathcal F}}_m(\varphi_t)\,\right\vert_{t=0}
\,=&\,\int_M\vert\nabla^m\nu\vert^2\,d\left(\left.\frac{\partial\mu_t}{\partial
    t}\,\right\vert_{t=0}\right) + \int_M\left.\frac{\partial}{\partial t}\vert
    \nabla^m\nu\vert^2\,\right\vert_{t=0}\,d\mu\\
\,=&\,\int_M\vert\nabla^m\nu\vert^2\, \HHH
\langle\nu\,\vert\,X\rangle\,d\mu\\
\,&\,+ \int_M\left.\frac{\partial}{\partial t}\left(g^{i_1j_1}\dots
    g^{i_mj_m}\nabla_{i_1\dots i_m}\nu\nabla_{j_1\dots
      j_m}\nu\right)\,\right\vert_{t=0}\,d\mu\,.
\end{align*}
Then, we need to compute the derivatives in the last term.\\
For the metric tensor $g_{ij}$ we have
\begin{align*}
\frac{\partial}{\partial t}g_{ij}\,=&\,\frac{\partial}{\partial
  t}\left\langle\frac{\partial\varphi}{\partial
    x_i}\,\right\vert\,\left.\frac{\partial\varphi}{\partial
    x_j}\right\rangle\\
\,=&\,\left\langle\frac{\partial X}{\partial
    x_i}\,\right\vert\,\left.\frac{\partial\varphi}{\partial
    x_j}\right\rangle+\left\langle\frac{\partial X}{\partial
    x_j}\,\right\vert\,\left.\frac{\partial\varphi}{\partial
    x_i}\right\rangle\\
\,=&\,\frac{\partial}{\partial x_i}\left\langle
  X\,\left\vert\,\frac{\partial\varphi}{\partial
    x_j}\right\rangle\right.
+\frac{\partial}{\partial x_j}\left\langle
  X\,\left\vert\,\frac{\partial\varphi}{\partial
    x_i}\right\rangle\right.
-2\left\langle
  X\,\left\vert\,\frac{\partial^2\varphi}{{\partial
    x_i}{\partial x_j}}\right\rangle\right.\\
\,=&\,a_{ij}(X)\,.
\end{align*}

Differentiating the formula $g_{is}g^{sj}=\delta_i^j$ we get
$$
\frac{\partial}{\partial  t}g^{ij}= - g^{is}\frac{\partial}{\partial
  t}g_{sl}g^{lj} = - g^{is}a_{sl}(X)g^{lj}\,.
$$

The derivative of the normal $\nu$ is given by
\begin{align*}
\frac{\partial}{\partial t}\nu\,=&\,\left\langle\frac{\partial\nu}{\partial
  t}\,\left\vert\,\frac{\partial\varphi}{\partial
    x_i}\right\rangle\right.\frac{\partial\varphi}{\partial
    x_j}g^{ij}\,=\,-\left\langle\nu\,\left\vert\,\frac{\partial^2\varphi}{\partial
    t\partial x_i}\right\rangle\right.\frac{\partial\varphi}{\partial
    x_j}g^{ij}\\
\,=&\,-\left\langle\nu\,\left\vert\,\frac{\partial X}{\partial
      x_i}\right\rangle\right.\frac{\partial\varphi}{\partial
  x_j}g^{ij}\,=\,-\nabla\left\langle\nu\,\left\vert\,X\right\rangle\right. +
\left.\left\langle\frac{\partial \nu}{\partial
      x_i}\,\right\vert\,X\right\rangle
\frac{\partial\varphi}{\partial  x_j}g^{ij}\\
\,=&\,-\nabla\left\langle\nu\,\left\vert\,X\right\rangle\right. +
\nabla\nu^\alpha X^\alpha=b(X)\,.
\end{align*}

Finally the derivative of the Christoffel symbols is
\begin{align*}
\frac{\partial}{\partial t}\Gamma_{jk}^i
=&\,\frac{1}{2}g^{il}\left\{
\frac{\partial}{\partial x_j}\left(\frac{\partial}{\partial t}g_{kl}\right) + 
\frac{\partial}{\partial x_k}\left(\frac{\partial}{\partial t}g_{jl}\right) - 
\frac{\partial}{\partial x_l}\left(\frac{\partial}{\partial
    t}g_{jk}\right)\right\}\\
&\,+\frac{1}{2}\frac{\partial}{\partial t}g^{il}\left\{
\frac{\partial}{\partial x_j}g_{kl} + 
\frac{\partial}{\partial x_k}g_{jl} - 
\frac{\partial}{\partial x_l}g_{jk}\right\}\\
=&\,\frac{1}{2}g^{il}\left\{
\nabla_j\left(\frac{\partial}{\partial t}g_{kl}\right) + 
\nabla_k\left(\frac{\partial}{\partial t}g_{jl}\right) - 
\nabla_l\left(\frac{\partial}{\partial t}g_{jk}\right)\right\}\\
&\,+\frac{1}{2}g^{il}\left\{
\frac{\partial}{\partial t}g_{kz}\Gamma_{jl}^z 
+\frac{\partial}{\partial t}g_{lz}\Gamma_{jk}^z
+\frac{\partial}{\partial t}g_{jz}\Gamma_{kl}^z 
+\frac{\partial}{\partial t}g_{lz}\Gamma_{jk}^z
-\frac{\partial}{\partial t}g_{jz}\Gamma_{kl}^z 
-\frac{\partial}{\partial t}g_{kz}\Gamma_{jl}^z\right\}\\
&\,-\frac{1}{2}g^{is}\frac{\partial}{\partial t}g_{sz}g^{zl}\left\{
\frac{\partial}{\partial x_j}g_{kl} + 
\frac{\partial}{\partial x_k}g_{jl} - 
\frac{\partial}{\partial x_l}g_{jk}\right\}\\
=&\,\frac{1}{2}g^{il}\left\{
\nabla_j\left(\frac{\partial}{\partial t}g_{kl}\right) + 
\nabla_k\left(\frac{\partial}{\partial t}g_{jl}\right) - 
\nabla_l\left(\frac{\partial}{\partial t}g_{jk}\right)\right\}\\
&\,+ g^{il}\frac{\partial}{\partial t}g_{lz}\Gamma_{jk}^z 
- g^{is}\frac{\partial}{\partial t}g_{sz}\Gamma_{jk}^z\\
=&\,\frac{1}{2}g^{il}\left\{
\nabla_j\left(\frac{\partial}{\partial t}g_{kl}\right) + 
\nabla_k\left(\frac{\partial}{\partial t}g_{jl}\right) - 
\nabla_l\left(\frac{\partial}{\partial t}g_{jk}\right)\right\}\\
=&\,\frac{1}{2}g^{il}\left\{
\nabla_ja_{kl}(X) + 
\nabla_ka_{jl}(X) - 
\nabla_la_{jk}(X)\right\}\,.
\end{align*}
Notice that all these derivatives are linear in the field $X$, since
the $a_{ij}(X)$ and $b(X)$ are such.

\begin{lemma}\label{eulerlemma}
If $a(X)=\frac{\partial}{\partial t}g$ is the
tensor defined before, for every covariant tensor $T=T_{i_1\dots i_l}$ 
we have
$$
\frac{\partial}{\partial t}\nabla^sT=\nabla^s\frac{\partial T}{\partial
  t}+\pol_{s-1}(T, \nabla a(X))
$$
where the constants in the polynomials $\pol_{s-1}(T, \nabla a(X))$
are universal.\\
Moreover, if the tensor $T$ is a function $f:M\to\R^k$ the last term
$\pol_{s-1}(f, \nabla a(X))$ can be substituted with another
polynomial $\widetilde{\pol}_{s-2}(\nabla f, \nabla a(X))$.
\end{lemma}

\begin{proof}
We prove the lemma by induction on $s\geq1$.\\
If $s=1$ then 
\begin{align*}
\frac{\partial}{\partial t}\nabla_{j}T_{i_1\dots i_l}
\,=&\,\frac{\partial}{\partial t}\left(\frac{\partial}{\partial
    x_{j}}T_{i_1\dots i_l}-\Gamma_{ji_z}^r T_{i_1\dots
    i_{z-1}ri_{z+1}\dots i_l}\right)\\
\,=&\,\frac{\partial}{\partial x_j}\frac{\partial}{\partial
    t}T_{i_1\dots i_l}-\Gamma_{ji_z}^r \frac{\partial}{\partial
    t}T_{i_1\dots i_{z-1}ri_{z+1}\dots i_l}\\
\,&\,-\frac{\partial}{\partial t}\Gamma_{ji_z}^rT_{i_1\dots
    i_{z-1}ri_{z+1}\dots i_l}\\
\,=&\,\nabla\frac{\partial T}{\partial
  t} + T * \nabla a(X)
\end{align*}
by the previous computation, hence
$$
\frac{\partial}{\partial t}\nabla T = \nabla\frac{\partial T}{\partial
  t}+\pol_0(T, \nabla a(X))
$$
and the initial case is proved.\\
Supposing the lemma holds for $s-1$, we have
\begin{align*}
\frac{\partial}{\partial t}\nabla^sT
\,=&\,\frac{\partial}{\partial t}\nabla(\nabla^{s-1}T)\\
\,=&\,\nabla\left(\frac{\partial}{\partial
    t}\nabla^{s-1}T\right) + \pol_0(\nabla^{s-1}T, \nabla a(X))\\
\,=&\,\nabla\left(\nabla^{s-1}\frac{\partial T}{\partial
    t} + \pol_{s-2}(T, \nabla a(X))\right)\\
\,&\,+\pol_0(\nabla^{s-1}T, \nabla a(X))\\
\,=&\,\nabla^{s}\frac{\partial T}{\partial t} + 
\nabla\pol_{s-2}(T, \nabla a(X))\\
\,&\,+\pol_0(\nabla^{s-1}T, \nabla a(X))\\
\,=&\,\nabla^{s}\frac{\partial T}{\partial t}+\pol_{s-1}(T, \nabla a(X))\\
\end{align*}
where we set
$$
\pol_{s-1}(T, \nabla a(X))=\nabla\pol_{s-2}(T, \nabla a(X)) +
\pol_0(\nabla^{s-1}T, \nabla a(X))\,.
$$
By this last formula, it is clear that the constants involved are
universal. Moreover, if $T$ is a function
$f:M\to\R^k$ then the term $\pol_0(f, \nabla a(X))$ vanishes and 
the same formula says that $\pol_{s-1}(f, \nabla a(X))$ does not
contain $f$ without being differentiated.
\end{proof} 

\begin{rem} In the following we will omit to underline that all the 
  coefficients of the polynomials $\pol_s$ and $\qol^s$ which will
  appear are {\em algebraic}, that is, they are the result
  of formal manipulations. In particular, such coefficients are
  independent of the manifold $(M,g)$ where the tensors are defined. 
  This is crucial in view of the geometry--independent estimates we
  want to obtain.
\end{rem}

\begin{prop}
The derivative
$$
\left.\frac{\partial}{\partial t}\left(g^{i_1j_1}\dots
    g^{i_mj_m}\nabla_{i_1\dots i_m}\nu\nabla_{j_1\dots
      j_m}\nu\right)\,\right\vert_{t=0}
$$
depends only on the vector field
$X=\left.\frac{\partial}{\partial t}\varphi_t\,\right\vert_{t=0}$ 
and such dependence is linear.\\
The first variation of ${\mathcal F}_m$
$$
\delta{\mathcal F}_m(\varphi)(\mathcal
I)=\left.\frac{d}{dt}{{\mathcal F}}_m(\varphi_t)\,\right\vert_{t=0}
$$
is a linear function of the field $X$.
\end{prop}

\begin{proof}
Distributing the derivative in $t$ on the terms of the product, we
have seen that the derivatives of the metric coefficients depends
linearly on $X$, it lasts to check the derivative of $\nabla_{i_1\dots 
i_m}\nu$.\\
By the last assertion of Lemma~\ref{eulerlemma}, we have
$$
\frac{\partial}{\partial t}\nabla^m\nu=\nabla^m\frac{\partial \nu}{\partial
  t}+\pol_{m-2}(\nabla \nu, \nabla a(X))
$$
and since $\frac{\partial \nu}{\partial
  t}=b(X)$ we get
$$
\frac{\partial}{\partial
  t}\nabla^m\nu=\nabla^m b(X) +\pol_{m-2}(\nabla \nu, \nabla a(X))
$$
which proves the first part of the lemma as $a(X)$ and $b(X)$ are
linear in $X$.\\
The second statement clearly follows by the previous computations and
the first part of the lemma.
\end{proof}

By this result, we can write $\delta{\mathcal F}_m(\varphi)(\mathcal
I)=\delta{\mathcal F}_m(\varphi)(X)$.
Now we want to prove that actually the first variation depends only 
on the normal component of the field  $X$, that is, 
$\langle\nu\,\vert\,X\rangle$, by linearity, it is clearly
sufficient to show that $\delta{\mathcal F}_m(\varphi)(X)=0$ for
every tangent vector field $X$. By the previous proposition, in order
to compute the derivative~\eqref{derivaz} we can choose any family
${\mathcal I}$ of immersions with $X=\left.\frac{\partial}{\partial
    t}\varphi_t\,\right\vert_{t=0}$.

Given a vector field $X$ along $M$ as a submanifold of $\Ri$ which is
tangent, there exists a tangent vector field $Y$ on $M$ such that
$d\varphi_p(Y(p))=X(p)$ for every $p\in M$.\\
Then we consider the smooth flow
$L(p,t):M\times(-\varepsilon,\varepsilon)\to M$ generated by $Y$ on
$M$ as the solution of the ODE's system
$$
\begin{cases}
\frac{\partial}{\partial t}L(p,t)\,=\,Y(L(p,t))\,,\\
L(p,0)\,=\,p
\end{cases}
$$
for every $p\in M$ and $t\in(-\varepsilon,\varepsilon)$, and we define
$\varphi_t(p)=\varphi(L(p,t))$.\\
Clearly $\varphi_0=\varphi$ and 
$$
\left.\frac{\partial}{\partial
    t}\varphi_t(p)\right\vert_{t=0}=
\left.d\varphi_{L(p,t)}\left(\frac{\partial}{\partial
      t}L(p,t)\right)\right\vert_{t=0}=d\varphi_p(Y(p))=X(p)\,,
$$
hence, using the family ${\mathcal I}=\{\varphi_t\}$ we have
$$
\delta{\mathcal F}_m(\varphi)(X)=\left.\frac{d}{dt}{{\mathcal
      F}}_m(\varphi_t)\,\right\vert_{t=0}\,.
$$
If $g_t$ is the metric tensor on $M$ induced by $\Ri$ via the
immersion $\varphi_t$, then the Riemannian manifolds $(M,g_t)$ and
$(M,g)$ are isometric for every $t\in(-\varepsilon,\varepsilon)$,
being $I(\cdot\,,t)=\varphi^{-1}\comp\varphi_t:(M,g_t)\to(M,g)$ an
isometry between them. Since  the functional ${\mathcal F}_m$ is
invariant by isometry, ${\mathcal F}_m(\varphi_t)$ does not depend on
$t$ and its derivative is zero.\\
By the previous discussion we have then the following proposition.

\begin{prop} The first variation $\delta{\mathcal F}_m(\varphi)(X)$
  depends only on $\langle\nu\,\vert\,X\rangle$.
\end{prop}

This means that we can suppose that $X$ is a normal field in studying 
$\delta{\mathcal F}_m(\varphi)(X)$, hence we can strengthen the
previous computations as follows,

\begin{align*}
\frac{\partial}{\partial t}g_{ij}\,
\,=&\, a_{ij}(X)\,=\,-2\left\langle
  X\,\left\vert\,\frac{\partial^2\varphi}{{\partial
    x_i}{\partial x_j}}\right\rangle\right.\,=\,
2\,\langle\nu\,\vert\,X\rangle h_{ij}\\
\frac{\partial}{\partial t}g^{ij}\,
\,=&\, - g^{is}\frac{\partial}{\partial
  t}g_{sl}g^{lj}\,=\, - 2\,\langle\nu\,\vert\,X\rangle h^{ij}\\
\frac{\partial}{\partial t}\nu\phantom{^{ij}}\,
\,=&\,-\nabla\left\langle\nu\left\vert\,X\right\rangle\right.\\
\frac{\partial}{\partial t}\Gamma_{jk}^i
\,=&\,g^{il}\left\{
\nabla_j(\langle\nu\,\vert\,X\rangle h_{kl}) + 
\nabla_k(\langle\nu\,\vert\,X\rangle h_{jl}) - 
\nabla_l(\langle\nu\,\vert\,X\rangle h_{jk})\right\}\\
\phantom{\frac{\partial}{\partial t}\Gamma_{jk}^i}
\,=&\,{\nabla\AAA} * \langle\nu\,\vert\,X\rangle
+\AAA * \nabla\langle\nu\,\vert\,X\rangle\,.
\end{align*}

Supposing $X$ normal, we have immediately the following modification
of Lemma~\ref{eulerlemma} substituting the tensor $a_{ij}(X)$ with 
$2\,\langle\nu\,\vert\,X\rangle h_{ij}$.

\begin{lemma}\label{eulerkey}
For every covariant tensor $T=T_{i_1\dots i_l}$, we have
$$
\frac{\partial}{\partial t}\nabla^sT=\nabla^s\frac{\partial T}{\partial
  t}+\pol_s(T, \AAA, \langle\nu\,\vert\,X\rangle)
$$
where in $\pol_s(T, \AAA, \langle\nu\,\vert\,X\rangle)$ the derivative 
$\nabla^s T$ does not appear.
If $T$ is a function $f:M\to\R^k$
$$
\frac{\partial}{\partial t}\nabla^sf=\nabla^s\frac{\partial f}{\partial
  t}+\pol_{s-1}(\nabla f, \AAA, \langle\nu\,\vert\,X\rangle)
$$
and $\pol_{s-1}(\nabla f, \AAA, \langle\nu\,\vert\,X\rangle)$ does not 
contain $\nabla^s f$.
\end{lemma}

This lemma and the fact that $\frac{\partial\nu}{\partial
t}=-\nabla\langle\nu\,\vert\,X\rangle$ lead to the following
proposition.

\begin{prop}
Letting $\{e_\alpha\}$ the canonical basis of $\Ri$ and setting
$\nu=\nu^\alpha e_\alpha\in\Ri$, we have
$$
\frac{\partial}{\partial t}\nabla_{i_1\dots i_m}\nu^\alpha=
-\nabla_{i_1\dots i_m}\nabla^\alpha\langle\nu\,\vert\,X\rangle 
+\pol_{m-1}(\nabla\nu, \AAA, \langle\nu\,\vert\,X\rangle)
$$
where we denoted with $\nabla^\alpha\langle\nu\,\vert\,X\rangle$ the
$\alpha$ component of the gradient $\nabla\langle\nu\,\vert\,X\rangle$
in the canonical basis of $\Ri$. Moreover, the derivative $\nabla^m\nu$ is
not present in $\pol_{m-1}(\nabla\nu, \AAA,
\langle\nu\,\vert\,X\rangle)$.
\end{prop}

We are finally ready to compute
\begin{align*}
\left.\frac{d}{dt}\int_M1+\vert
    \nabla^m\nu\vert^2\,d\mu_t\right\vert_{t=0}
\,=&\,\int_M\left(1+\vert
    \nabla^m\nu\vert^2\right)\HHH \langle\nu\,\vert\,X\rangle\,d\mu\\
&\,+\int_M g^{i_1j_1}\dots\frac{\partial}{\partial t}g^{i_kj_k}\dots
g^{i_mj_m}\nabla_{i_1\dots i_m}\nu \nabla_{j_1\dots j_m}\nu\,d\mu\\
&\,-2\int_M g^{i_1j_1}\dots g^{i_mj_m}\, 
\nabla_{i_1\dots i_m}\nabla^\alpha\langle\nu\,\vert\,X\rangle\,
\nabla_{j_1\dots  j_m}\nu^\alpha\,d\mu\\
&\,+2\int_M \nabla^m\nu * \pol_{m-1}(\nabla\nu, \AAA,
\langle\nu\,\vert\,X\rangle)\,d\mu\\
\,=&\,\int_M\left(1+\vert
    \nabla^m\nu\vert^2\right)\HHH \langle\nu\,\vert\,X\rangle\,d\mu\\
&\,+2m\int_M \nabla^m\nu * \nabla^m\nu *
\AAA \langle\nu\,\vert\,X\rangle\,d\mu\\
&\,-2\int_M g^{i_1j_1}\dots g^{i_mj_m}\, 
\nabla_{i_1\dots i_m}\nabla^\alpha\langle\nu\,\vert\,X\rangle\,
\nabla_{j_1\dots  j_m}\nu^\alpha\,d\mu\\
&\,+\int_M \pol_{m-1}(\nabla^m\nu, \nabla\nu, \AAA,
\langle\nu\,\vert\,X\rangle)\,d\mu\,.
\end{align*}
Now, in order to ``carry away'' derivatives from
$\langle\nu\,\vert\,X\rangle$ in the last integral, we integrate by
parts with the divergence theorem, ``moving'' all the derivatives on
the other terms of the products. Hence, we can rewrite it as
$$
\int_M\pol_{2m-2}(\nabla\nu, \nabla\nu, \AAA) 
\langle\nu\,\vert\,X\rangle\,d\mu\,,
$$
which is equal to 
$$
\int_M\qol^{2m+1}(\nabla\nu, \AAA) 
\langle\nu\,\vert\,X\rangle\,d\mu
$$
with the conventions of Section~\ref{prelim}.\\
Since also the second integral has this form, collecting them
together, we obtain
\begin{align*}
\left.\frac{d}{dt}\int_M1+\vert
    \nabla^m\nu\vert^2\,d\mu_t\right\vert_{t=0}
\,=&\,\int_M\HHH \langle\nu\,\vert\,X\rangle\,d\mu
+\int_M\qol^{2m+1}(\nabla\nu, \AAA) \langle\nu\,\vert\,X\rangle\,d\mu\\
&\,-2\int_M g^{i_1j_1}\dots g^{i_mj_m}\, 
\nabla_{i_1\dots i_m}\nabla^\alpha\langle\nu\,\vert\,X\rangle\,
\nabla_{j_1\dots  j_m}\nu^\alpha\,d\mu\,.
\end{align*}
Finally, we deal with this last term. First, 
by the divergence theorem it can be transformed in
$$
-2(-1)^{m}\int_M \nabla^\alpha\langle\nu\,\vert\,X\rangle\,
\nabla^{j_m\dots j_1}\nabla_{j_1\dots
  j_m}\nu^\alpha\,d\mu\,,
$$
second, using the tangential divergence formula~\eqref{divtang},
it is equal to 
$$
2(-1)^{m}\int_M \langle\nu\,\vert\,X\rangle\,
\nabla^\alpha\nabla^{j_m\dots j_1}\nabla_{j_1\dots
  j_m}\nu^\alpha\,d\mu + \int_M \qol^{2m+1}(\nabla\nu, \AAA)
\langle\nu\,\vert\,X\rangle\,d\mu\,,
$$
where the extra term $\qol^{2m+1}(\nabla\nu, \AAA)
\langle\nu\,\vert\,X\rangle$, which has a differentiation order lower
than the first term, comes from the product with the mean curvature
in the tangential divergence formula.\\
Notice now that the permutation of derivatives introduces
additional lower order terms of the form
$$
\int_M \qol^{2m+1}(\nabla\nu, \AAA) \langle\nu\,\vert\,X\rangle\,d\mu
$$
by formulas~\eqref{ichange}, hence we get
$$
2(-1)^{m}\int_M \langle\nu\,\vert\,X\rangle\, 
\nabla^{j_1}\nabla_{j_1}\dots\nabla^{j_m}\nabla_{j_m}\nabla^\alpha\nu^\alpha\,d\mu
+ \int_M \qol^{2m+1}(\nabla\nu, \AAA)
\langle\nu\,\vert\,X\rangle\,d\mu
$$
that is,
$$
2(-1)^{m}\int_M \langle\nu\,\vert\,X\rangle \overset{\text{$m$~{\rm
      times}}}{\overbrace{\Delta\Delta\dots\Delta}}
\nabla^\alpha\nu^\alpha\,d\mu + \int_M \qol^{2m+1}(\nabla\nu, \AAA)
\langle\nu\,\vert\,X\rangle\,d\mu\,.
$$
By Gauss--Weingarten relations~\eqref{gwein}, we have
$$
\nabla^\alpha\nu^\alpha=\frac{\partial\varphi^\alpha}{\partial
  x_i}g^{ij}h_{jl}g^{ls}\frac{\partial\varphi^\alpha}{\partial x_s}=
g^{ij}h_{jl}g^{ls}g_{si}=g^{ij}h_{ji}= \HHH\,,
$$
so we conclude
\begin{align*}
\delta{\mathcal F}_m(\varphi)(X)
\,=&\,\int_M\HHH \langle\nu\,\vert\,X\rangle\,d\mu
+\int_M\qol^{2m+1}(\nabla\nu, \AAA) \langle\nu\,\vert\,X\rangle\,d\mu\\
&\,+2(-1)^{m}\int_M \overset{\text{$m$~{\rm times}}}
{\overbrace{\Delta\Delta\dots\Delta}} \HHH 
\langle\nu\,\vert\,X\rangle\,d\mu\\
\,=&\,\int_M \qol^1(\AAA) 
\langle\nu\,\vert\,X\rangle\,d\mu
+\int_M\qol^{2m+1}(\nabla\nu, \AAA) 
\langle\nu\,\vert\,X\rangle\,d\mu\\
&\,+2(-1)^{m}\int_M \overset{\text{$m$~{\rm times}}}
{\overbrace{\Delta\Delta\dots\Delta}} \HHH 
\langle\nu\,\vert\,X\rangle\,d\mu\,.
\end{align*}

By the previous discussion this formula holds in general for every
vector field $X$ along $M$. We summarize all these facts in the
following theorem.

\begin{teo}\label{fine} For any $m\geq1$ the first variation of the
  functional ${\mathcal F}_m$ is given by
$$
\delta{\mathcal F}_m(\varphi)(X)=\int_M \EEE_m(\varphi) 
\langle\nu\,\vert\,X\rangle\,d\mu
$$
where the function $\EEE_m(\varphi)$ has the form
\begin{equation*}
\EEE_m(\varphi)=2(-1)^{m}\overset{\text{$m$~{\rm times}}}
{\overbrace{\Delta\Delta\dots\Delta}} \HHH + \qol^{2m+1}(\nabla\nu,
\AAA) + \qol^1(\AAA)\,.
\end{equation*}
\end{teo}

\section{Gradient Flow and Small Time Existence}

Suppose that $\varphi_0:M\to\R^{n+1}$ is smooth immersion of an
$n$--dimensional hypersurface $M$ which is compact, connected 
and has empty boundary.

We look for a smooth function $\varphi:M\times[0,T)$ such that
\begin{enumerate}
\item the map $\varphi_t=\varphi(\cdot,t):M\to\R^{n+1}$ is an
  immersion;
\item the following partial differential equation is satisfied
\begin{equation*}
\frac{\partial\varphi}{\partial t}(p,t)=-\EEE_m(\varphi_t)(p)\nu(p,t)\,.
\end{equation*}
\end{enumerate}
If we have a solution, then we say that the hypersurfaces
$M_t=(M,g_t)$, where $g_t$ is the induced metric on
  $M$, evolve by the gradient flow of the functional ${\mathcal
    F}_m$.

The small time existence of such flow is a slight
modification of the following result of Polden (see~\cite{polden2},
Thm.~2.5.2, Sec.~2 or~\cite{huiskpold}).

\begin{teo}\label{smalltime}
For any smooth hypersurface immersion $\varphi_0:M\to N$, with $N$ a
smooth $(n+1)$--dimensional Riemannian manifold, there
exists a unique solution to the flow problem
$$
\frac{\partial \varphi}{\partial
  t}=\biggl((-1)^{s+1}\overset{\text{$s$~{\rm times}}}
{\overbrace{\Delta\Delta\dots\Delta}} \HHH + \Phi(\varphi, \nu, \AAA,
\nabla \AAA, \,\dots\,, \nabla^{2s-1}\AAA)\biggr)\nu
$$
defined on some interval $0\leq t<T$ and taking $\varphi_0$ as its initial
value.
\end{teo}

Looking at Polden's proof, it is possible to allow the function 
$\Phi$ to depend also on the metric $g$, moreover the 
covariant derivatives of the normal $\nu$, using induction and the
Gauss--Weingarten relations~\eqref{gwein}, can be expressed in terms of
the covariant derivatives of the curvature (see the proof of
Lemma~\ref{covacurv}).\\
Hence, we can conclude that there exists a small time solution of the
problem
$$
\frac{\partial \varphi}{\partial
  t}=\biggl((-1)^{m+1}\overset{\text{$m$~{\rm times}}}
{\overbrace{\Delta\Delta\dots\Delta}}\,\,H\,
+\Phi(\varphi, g, \AAA, \nu, \nabla\AAA, \nabla\nu, \,\dots\,,
\nabla^{2m-1}\AAA, \nabla^{2m}\nu)\biggr)\nu
$$
which includes our case up to a constant multiplying the leading
term. Since such a constant can be eliminated by a time--only
rescaling and since a smooth evolution of an immersed compact manifold 
clearly remains an immersion at least for some positive time, 
we have a small time existence and uniqueness result for the gradient
flow of ${\mathcal F}_m$ with every initial hypersurface.

\section{A Priori Estimates}

To prove long time existence we need a priori estimates
on the second fundamental form and its derivatives which are 
obtained via Sobolev and Gagliardo--Nirenberg interpolation
inequalities for functions defined on $M_t$.\\
Since the hypersurfaces are moving, also the constants appearing in
such inequalities change during the flow, hence, before proceeding with
the estimates, we need some uniform control on them.

In this section we see that if the integer $m$ larger than 
$\left[\frac{n}{2}\right]$ then we have a uniform control, independent
of time,  on the $L^{n+1}$ norm of the second fundamental form; this
is a crucial point where such hypothesis is necessary. 
This fact will allow us to show in the next section that also the above
constants are uniformly bounded during the flow.\\
In the last part of the section, using an inequality of Michael and
Simon, we prove also an a priori lower bound on the volume of the
evolving hypersurfaces.

By the very definition of the flow, the value of the functional
${\mathcal F}_m$ decreases in time, since
$$
\frac{d}{dt}{\mathcal F}_m(\varphi_t)=
-\int_{M}\left[\EEE_m(\varphi_t)\right]^2\,d\mu_t\leq0\,,
$$
hence, as long as the flow remains smooth, we have the 
uniform estimate
\begin{equation}\label{pippo}
\int_{M}1+\vert\nabla^m\AAA\vert^2\,d\mu_t={\mathcal
  F}_m(\varphi_t)\leq{\mathcal F}_m(\varphi_0)
\end{equation}
for every $t\geq0$.

Now we want to prove that if $m>\left[\frac{n}{2}\right]$, this 
estimate implies that the $L^{n+1}(\mu_t)$ norms of the second
fundamental form $\AAA$ of $M_t$ are uniformly bounded independently
of time.

Our starting point is the following universal interpolation type
inequalities for tensors.

\begin{prop}\label{univ1gen} Suppose that $(M,g)$ is a smooth and compact
  $n$--dimensional Riemannian manifold without boundary and $\mu$ the
  measure associated to $g$.\\
Then for every covariant tensor $T$ and exponents $q\in[1,+\infty)$
and $r\in[1,+\infty]$, we have
\begin{equation}\label{ccc3}
\Vert\nabla^j T\Vert_{L^{p}{(\mu)}}\leq
C\Vert\nabla^{s}T\Vert^{\frac{j}{s}}_{L^q{(\mu)}}\Vert
T\Vert_{L^r(\mu)}^{\frac{s-j}{s}}\qquad \forall j\in[0,s]\,,
\end{equation}
with
$$
\frac{1}{p}=\frac{j}{sq}+\frac{s-j}{sr}\,,
$$
where the constant $C$ depends only on $n$, $s$, $j$, $p$, $q$, $r$
and not on the metric or the geometry of $M$.
\end{prop}
The proof of the case $r=+\infty$ can be found in~\cite{hamilton1},
Sec.~12, along the same lines also the case $r<+\infty$ follows 
(see also~\cite{aubin0}, Chap.~3, Sec.~7.6).

Suppose that $M$ is orientable and that $g$ is the metric induced by 
the immersion $\varphi:M\to\Ri$, let $\nu$ be a global unit normal
vector field on $M$.\\
If in~\eqref{ccc3} we consider $T=\nu$, $s=m$,
$j=1$, $q=2$ and $r=+\infty$, then we have $\vert T\vert=1$ and
$p=2m$, hence
\begin{equation*}
\Vert\nabla\nu\Vert_{L^{2m}(\mu)}\leq
C\Vert\nabla^{m}\nu\Vert^{\frac{1}{m}}_{L^2(\mu)}\,,
\end{equation*}
for a constant $C=C(n, m)$.\\
Since by~\eqref{gwein} $\vert\nabla\nu\vert=\vert\AAA\vert$, we
conclude
$$
\int_M\vert{\AAA}\vert^{2m}\,d\mu\leq 
C \int_M\vert\nabla^{m}\nu\vert^2\,d\mu\leq 
C {\mathcal F}_m(\varphi)\,.
$$
If $M$ is not orientable, then there exists a two--fold Riemannian
covering $\widetilde{M}$ of $M$, with a locally isometric projection
map $\pi:\widetilde{M}\to M$ which is orientable and immersed in
$\Ri$ via the map $\varphi\comp\pi$. Repeating the previous 
argument for $\widetilde{M}$ we get
$$
\int_{\widetilde{M}}\vert{\AAA}\vert^{2m}\,d\widetilde{\mu}\leq C 
\int_{\widetilde{M}}\vert\nabla^{m}\nu\vert^2\,d\widetilde{\mu}\,.
$$
Since $\pi$ is a local isometry and noticing that the global unit normal
field on $\widetilde{M}$ gives locally a unit normal field on $M$,
all the quantities which appear inside the integrals above do not
change passing from $\widetilde{M}$ to $M$, only when we integrate we
need to take into account the two--fold structure of the covering. This
means that for every smooth function $u:M\to\R$ we have
$$
\int_{\widetilde{M}} u\comp\pi\,d\widetilde{\mu} = 2 \int_M u\,d\mu\,.
$$
Hence, we deduce
$$
2\int_M \vert{\AAA}\vert^{2m}\,d\mu\leq 
2 \, C \int_M \vert\nabla^{m}\nu\vert^2\,d\mu\leq 
2 \, C {\mathcal F}_m(\varphi)
$$
which clearly gives the same estimate as in the orientable case.

As $2m>2\left[\frac{n}{2}\right]\geq{n+1}$, we have
\begin{equation}\label{aequa15}
\int_M\vert{\AAA}\vert^{n+1}\,d\mu\leq \left(\int_M
  \vert\AAA\vert^{2m}\,d\mu\right)^{\frac{n+1}{2m}}\left(\text{{\rm
      {Vol}} $M$}\right)^{\frac{2m-n-1}{2m}}\leq C {\mathcal
  F}_m(\varphi)
\end{equation}
with a constant $C=C(n, m)$.

Finally we show that also the volume of $M$ is well controlled by
the value of ${\mathcal F}_m(\varphi)$ under the hypothesis
$m>\left[\frac{n}{2}\right]$.\\
The bound from above is obvious, the bound from below in dimension
$n>1$ can be obtained via the following universal Sobolev inequality
due to Michael and Simon (see~\cite{micsim,simon}).

\begin{prop}
Let $\varphi:M\to\Ri$ be an immersion of an
$n$--dimensional, compact hypersurface without boundary. On $M$ we
consider the Riemannian metric induced by $\Ri$ and the corresponding
measure $\mu$.\\ 
Then, there exists a constant $C=C(n, p)$ depending only on the dimension $n$
and the exponent $p$ such that, for every smooth function $u:M\to\R$
\begin{equation}\label{aequa16}
\left(\int_M\vert u\vert^{p^*}\,d\mu\right)^{1/p^*}
\leq C(n, p)\left(\int_M\vert\nabla u\vert^p\,d\mu +
\int_M \vert \HHH u\vert^p\,d\mu\right)^{1/p}\,,
\end{equation}
where $p\in[1,n)$, $n>1$ and $p^*=\frac{np}{n-p}$.
\end{prop}

Considering the function $u:M\to\R$ constantly equal to 1 in the 
inequality for $p=1$, and taking in account~\eqref{aequa15}, we get
\begin{align*}
\left(\text{{\rm {Vol}} $M$}\right)^{\frac{n-1}{n}}\leq &C\int_M \vert
{\HHH}\vert\,d\mu\\
\leq &C\Vert{\AAA}\Vert_{L^{n+1}{(\mu)}}\left(\text{{\rm {Vol}} 
$M$}\right)^{\frac{n}{n+1}}\\
\leq &C{\mathcal
  F}_m(\varphi)^{\frac{1}{n+1}}\left(\text{{\rm {Vol}} $M$}\right)^{\frac{n}{n+1}}\,.
\end{align*}
Dividing both members by $\left(\text{{\rm {Vol}} $M$}\right)^{\frac{n-1}{n}}$, as 
${\frac{n}{n+1}}>\frac{n-1}{n}$ we conclude
\begin{equation*}
1\leq C{\mathcal F}_m(\varphi)^{\frac{1}{n+1}}\left(\text{{\rm {Vol}} $M$}\right)^{\frac{1}{n(n+1)}}
\end{equation*}
that is,
\begin{equation*}
\frac{C}{{\mathcal F}_m(\varphi)^{n}}\leq\text{{\rm{Vol}} $M$}\leq{\mathcal F}_m(\varphi)
\end{equation*}
for a constant $C=C(n, m)$.
\begin{rem} With the same argument, it follows that 
  also $\Vert{\AAA}\Vert_{L^{n+1}{(\mu)}}$ 
  can be controlled above and below with  ${\mathcal F}_m(\varphi)$ 
  and that the functional ${\mathcal F}_m$ is
  uniformly bounded from below by a constant greater than zero.
\end{rem}

In the special case $n=1$, we recall that for every closed curve
$\gamma:\SS^1\to\R^2$ in the plane the integral of the modulus of its
curvature is at least $2\pi$, then
$$
2\pi\leq\int_{\SS^1}\vert{\AAA}\vert\,ds\leq\left(\int_{\SS^1}\vert{\AAA}
\vert^2\,ds\right)^{1/2}\sqrt{\text{${\mathrm {Length}}$ $\gamma$}}\leq
C \sqrt{{\mathcal F}_m(\gamma)}\sqrt{\text{${\mathrm {Length}}$
    $\gamma$}}\,.
$$
Hence,
$$
\frac{C}{{\mathcal F}_m(\gamma)}\leq{\text{${\mathrm {Length}}$
    $\gamma$}}\leq{{\mathcal F}_m(\gamma)}
$$
with $C=C(m)$.

Putting together all these inequalities and the uniform
estimate~\eqref{pippo} we obtain the following result.

\begin{prop}\label{n+1} As long as the flow by the gradient of
  ${\mathcal F}_m$ of a hypersurface in $\Ri$ exists, we have the
  estimates
\begin{equation*}
\Vert{\AAA}\Vert_{L^{n+1}(\mu_t)}\leq C_1 <+\infty
\end{equation*}
\begin{equation*}
0 < C_2 \leq \text{{\rm {Vol}} $M_t$}\leq C_3 < +\infty
\end{equation*}
where the three constants $C_1$, $C_2$ and  $C_3$ are independent of 
time.\\
They depend only on  $n$, $m$ and the value of ${\mathcal F}_m$ for
the initial hypersurface.
\end{prop}

\section{Interpolation Inequalities for Tensors}\label{interptensor}

As we said, we show now that the uniform bound on the $L^{n+1}$ norm
of the second fundamental form implies that the constants involved in
some Sobolev and Gagliardo--Nirenberg interpolation type
inequalities are also equibounded.

Recalling inequality~\eqref{aequa16}, we have
\begin{equation}\label{ssss1}
\Vert u\Vert_{L^{p^*}(\mu)}\leq C(n, p)\left(\Vert\nabla
  u\Vert_{L^p(\mu)} +\Vert \HHH u\Vert_{L^p(\mu)}\right)
\end{equation}
for every $u\in C^1(M)$, where ${p^*}=\frac{np}{n-p}$ and $p\in[1,n)$.

\begin{prop}\label{sobsim1000} If  the manifold $(M,g)$ satisfies
  {\rm {Vol}} $ M +\Vert\HHH\Vert_{L^{n+\delta}(\mu)}\leq B$ for some
  $\delta>0$ then for every $p\in[1,n)$,
\begin{equation*}
\Vert u\Vert_{L^{p^*}(\mu)}\leq 
C\left(\Vert\nabla u\Vert_{L^{p}(\mu)} 
+ \Vert u\Vert_{L^p(\mu)}\right)\qquad\forall u\in C^1(M)\,,
\end{equation*}
where the constant $C$ depends only on $n$, $p$, $\delta$ and $B$.
\end{prop}

\begin{proof}
Applying H\"older inequality to the last term of 
inequality~\eqref{ssss1}, we get
\begin{equation*}
\Vert u\Vert_{L^{p^*}(\mu)}\leq C(n, p)\Vert\nabla
  u\Vert_{L^p(\mu)} + C(n, p, \delta, B)\Vert
  u\Vert_{L^{\widetilde{p}}(\mu)}
\end{equation*}
where $\widetilde{p}$ is given by
$$
\widetilde{p}=\frac{p(n+\delta)}{n+\delta-p}=p^*\frac{n(n+\delta)}
{n(n+\delta)+   p^*\delta}\,,
$$
then $p<\widetilde{p}<p^*$.\\
Hence, we can interpolate $\Vert u\Vert_{L^{\widetilde{p}}(\mu)}$ 
between a small fraction of $\Vert u\Vert_{L^{p^*}(\mu)}$ 
and a possibly large multiple of $\Vert u\Vert_{L^p(\mu)}$,
\begin{equation*}
\Vert u\Vert_{L^{p^*}(\mu)}\leq C(n, p)\Vert\nabla
  u\Vert_{L^p(\mu)} + C(n, p, \delta, B)
\left(\varepsilon\Vert u\Vert_{L^{p^*}(\mu)}+ 
C(\varepsilon, p)\Vert u\Vert_{L^p(\mu)}\right)\,.
\end{equation*}
Choosing $\varepsilon>0$ such that $\varepsilon C(n, p, \delta,
B)\leq1/2$ and collecting terms we obtain
\begin{equation*}
\Vert u\Vert_{L^{p^*}(\mu)}\leq C(n, p, \delta, B)\left(\Vert\nabla
  u\Vert_{L^p(\mu)} + \Vert u \Vert_{L^p(\mu)}\right)\,.
\end{equation*}
\end{proof}

When $p>n$ we prove the following $L^{\infty}$ result
(see also~\cite{kuschat1}, Thm.~5.6).

\begin{prop} If  the manifold $(M,g)$ satisfies
  {\rm {Vol}} $ M +\Vert\HHH\Vert_{L^{n+\delta}(\mu)}\leq B$ for some
  $\delta>0$ then for every $p>n$, we have
\begin{equation*}
\max_M\vert u\vert\leq 
C\left(\Vert\nabla u\Vert_{L^{p}(\mu)} 
+ \Vert u\Vert_{L^p(\mu)}\right)\qquad\forall u\in C^1(M)\,,
\end{equation*}
where the constant $C$ depends only on $n$, $p$, $\delta$ and $B$.
\end{prop}

\begin{proof}
Suppose first that $M$ is embedded and $n+\delta\geq p>n$, clearly $\Vert
\HHH\Vert_{L^p(\mu)}$ is bounded by a value depending on the constant 
$B$.\\
We consider $M$ as a subset of $\Ri$ via the embedding  $\varphi$ and 
$\mu$ as a measure on $\Ri$ which is supported on $M$. Then 
the following result holds (\cite{simon}, Thm.~17.7): let $B_\rho(x)$
be the ball of radius $\rho$ centered at $x$ in $\Ri$, for every
$0<\sigma<\rho<+\infty$ we have
$$
\left(\frac{\mu(B_\sigma(x))}{\sigma^n}\right)^{1/p}\leq
\left(\frac{\mu(B_\rho(x))}{\rho^n}\right)^{1/p} 
+ C(n, p, \delta, B)\left(\rho^{1-n/p}-\sigma^{1-n/p}\right)\,.
$$
Hence, 
$$
\left(\frac{\mu(B_\sigma(x))}{\sigma^n}\right)^{1/p}\leq
\frac{C_1}{\rho^{n/p}} + C_2\rho^{1-n/p}\,,
$$
and choosing $\rho=1$, for every $0<\sigma<1$ we get the inequality
$$
\mu(B_\sigma(x))\leq C(n, p, \delta, B)\sigma^n\,.
$$
Then we need the following formula which is proved in~\cite{simon},
Sec.~18, as a consequence of the tangential divergence 
formula~\eqref{divtang}.\\
For every $0<\sigma<\rho<+\infty$ we have
$$
\frac{\int_{B_\sigma(x)}u\,d\mu}{\sigma^n}\leq
\frac{\int_{B_\rho(x)}u\,d\mu}{\rho^n}+\int_\sigma^\rho\tau^{-n-1}
\int_{B_\tau(x)}r(\vert\nabla u\vert+\vert u\HHH\vert)\,d\mu(y)\,d\tau
$$
where $r=\vert x - y\vert$ and $u$ is any smooth non negative function.

Noticing that $r\leq\tau$ and using H\"older inequality we estimate
\begin{align*}
\frac{\int_{B_\sigma(x)}u\,d\mu}{\sigma^n}
\,\leq&\,\frac{\int_{B_\rho(x)}u\,d\mu}{\rho^n}+
\left(\int_M \vert\nabla
u\vert^p+\vert u\HHH\vert^p \,d\mu\right)^{1/p}
\int_\sigma^\rho\tau^{-n}\mu(B_\tau(x))^{1-1/p}\,d\tau\\
\,\leq&\,\int_{B_{1}(x)}u\,d\mu+
C \left(\Vert\nabla
u\Vert_{L^p(\mu)}+\Vert u\HHH\Vert_{L^p(\mu)}\right)
\int_\sigma^{1}\tau^{-n}\tau^{n-n/p}\,d\tau
\end{align*}
where in the last passage we set $\rho=1$ used the previous estimate on
$\mu(B_\tau(x))$. The function $\tau^{-n/p}$ is integrable since $p>n$
and we get
$$
\frac{\int_{B_\sigma(x)}u\,d\mu}{\sigma^n}
\leq {\int_{B_{1}(x)}u\,d\mu} +
C \left(\Vert\nabla
u\Vert_{L^p(\mu)}+\Vert u\HHH\Vert_{L^p(\mu)}\right)
\frac{1-\sigma^{1-n/p}}{1-n/p}\,,
$$
now sending $\sigma$ to zero, on the left side we obtain the value of
$u(x)$ times $\omega_n$ which is the volume of the unit ball of
$\R^n$, hence
\begin{align*}
\omega_n u(x)
\,\leq&\, {\int_{B_{1}(x)}u\,d\mu} +
C \left(\Vert\nabla
u\Vert_{L^p(\mu)}+\Vert u\HHH\Vert_{L^p(\mu)}\right)\\
\,\leq&\, C(n, p, \delta, B)\left(\Vert
  u\Vert_{L^1(\mu)}+ \Vert\nabla u\Vert_{L^p(\mu)}+\Vert
  u\HHH\Vert_{L^p(\mu)}\right)\,.
\end{align*}
For a general $u$ we apply this inequality to the function $u^2$, thus
\begin{align*}
u^2(x)\,&\leq\, C\left(\int_M \vert u\vert^2\,d\mu + 
\left(\int_M \vert u \nabla u\vert^p\,d\mu\right)^{1/p} +
\left(\int_M \vert u^2 \HHH \vert^{p}\,d\mu\right)^{1/p}\right)\\
\,&\leq\, C\max_M \vert u\vert \left(\int_M \vert u\vert\,d\mu + 
\left(\int_M \vert \nabla u\vert^p\,d\mu\right)^{1/p} +
\left(\int_M \vert u \HHH \vert^{p}\,d\mu\right)^{1/p}\right)\,.
\end{align*}
Since $x\in\Ri$ was arbitrary we conclude that
$$
\max_M\vert u\vert \leq 
C(n, p, \delta, B)\left(\Vert
  u\Vert_{L^1(\mu)}+ \Vert\nabla u\Vert_{L^p(\mu)}+\Vert
  u\HHH\Vert_{L^p(\mu)}\right)\,.
$$
for a constant $C$ depending on $n, p$, $\delta$ and $B$.\\
If $M$ is only immersed, we consider the embeddings of $M$ in
$\Ri\times\R^k$ given by the map
$\varphi\times\varepsilon\psi:M\to\Ri\times\R^k$, where
$\psi:M\to\R^k$ is an embedding of $M$ in some Euclidean space. Then,
repeating the previous argument (it is possible since the starting
inequalities from~\cite{simon} hold for embeddings in any $\R^l$) we
will get the same conclusion with a constant
$C_\varepsilon$. Finally, as $C_\varepsilon$ depends only on Vol $M$
and $\HHH$, and all the geometric quantities converge uniformly when
$\varepsilon$ goes to zero, we conclude that the  inequality holds
also in the immersed case.

Now, given any $p>n$, we choose $\widetilde{p}=\frac{1}{2}\min\{n+p,
2n+\delta\}$, then clearly $n<\widetilde{p}<\min\{p, n+\delta/2\}$. 
By the inequality above we have
$$
\max_M\vert u\vert\leq C(n, \widetilde{p}, \delta, B)
\left(\Vert u\Vert_{L^1(\mu)} + 
\Vert\nabla u\Vert_{L^{\widetilde{p}}(\mu)}+
\Vert u\HHH\Vert_{L^{\widetilde{p}}(\mu)}\right)\,,
$$
then using H\"older inequality and an interpolation argument as in the
proof of Proposition~\ref{sobsim1000} we get
$$
\max_M\vert u\vert\leq C(n, \widetilde{p}, \delta, B)
\left(\Vert u\Vert_{L^1(\mu)} + 
\Vert\nabla u\Vert_{L^{\widetilde{p}}(\mu)}+
\Vert u \Vert_{L^{{p}}(\mu)}\right)\,.
$$
Applying again H\"older inequality, as $\widetilde{p}<p$, we conclude
that
$$
\max_M\vert u\vert\leq C(n, \widetilde{p}, \delta, B)
\left(\Vert\nabla u\Vert_{L^{{p}}(\mu)} +
\Vert u \Vert_{L^{{p}}(\mu)}\right)\,,
$$
which gives the thesis since $\widetilde{p}$ depends only on $n$, $p$
and $\delta$.
\end{proof}

We now extend these propositions to tensors
(see~\cite{aubin0}, Prop.~2.11 and also~\cite{cantor1,cantor2}). 
Since $\vert T \vert$ is not necessarily smooth we apply the
previous inequalities first to the smooth functions $\sqrt{\vert
  T\vert^2 + \varepsilon^2}$, converging to $\vert T \vert$ when
$\varepsilon\to0$. As
$$
\left\vert \,\nabla \sqrt{\vert T\vert^2 + \varepsilon^2}\,\right\vert\,=\,
\left\vert \,\frac{\langle \nabla T , T\rangle}{\sqrt{\vert T\vert^2 +
      \varepsilon^2}}\, \right\vert
\,\leq\,\frac{\vert T \vert}{\sqrt{\vert T\vert^2 +
      \varepsilon^2}}\,\vert\nabla T\vert\,\leq \,\vert\nabla T\vert
$$
we get then easily the following result.

\begin{prop}\label{pnnp} If  the manifold $(M,g)$ satisfies
  {\rm {Vol}} $ M +\Vert\HHH\Vert_{L^{n+\delta}(\mu)}\leq B$ for some
  $\delta>0$   then for every covariant tensor $T=T_{i_1 \dots i_l}$ we
  have,
\begin{align}
\label{firstens}\Vert T\Vert_{L^{p^*}(\mu)}\leq\,&\, 
C\left(\Vert\nabla T\Vert_{L^{p}(\mu)} 
+ \Vert T\Vert_{L^p(\mu)}\right)\qquad\text{ if $1\leq p<n$,}\\
\label{sectens}\max_M\vert T\vert\leq\,&\, 
C\left(\Vert\nabla T\Vert_{L^{p}(\mu)} 
+ \Vert T\Vert_{L^p(\mu)}\right)\qquad\text{ if $p>n$,}
\end{align}
where the constants depend only on $n$, $l$, $p$, $\delta$ and
$B$.
\end{prop}

We define the Sobolev norm of a tensor $T$ on $(M,g)$ as
$$
\Vert T\Vert_{W^{s,q}(\mu)}=\sum_{i=0}^s \Vert\nabla^i
T\Vert_{L^q(\mu)}\,.
$$

\begin{cor} In the same hypothesis on $(M,g)$ we have
\begin{equation}\label{a=1}
\Vert\nabla^jT\Vert_{L^{p}(\mu)}\leq 
C\Vert T\Vert_{W^{s,q}(\mu)}\qquad \text{ with }\qquad
\frac{1}{p}=\frac{1}{q}-\frac{s-j}{n}>0\,,
\end{equation}
\begin{equation}\label{a=1inf}
\max_M\vert\nabla^jT\vert\leq 
C\Vert T\Vert_{W^{s,q}(\mu)}\qquad \text{ when }\qquad 
\frac{1}{q}-\frac{s-j}{n}<0\,.
\end{equation}
The constants depend only on $n$, $l$, $s$, $j$, $p$, $q$, $\delta$
and $B$.
\end{cor}

\begin{proof} By inequality~\eqref{firstens} applied to the tensor
  $\nabla^j T$ we get
\begin{align*}
\Vert \nabla^j T\Vert_{L^{p}(\mu)}
\,&\,\leq\,C\left(\Vert\nabla^{j+1} T\Vert_{L^{p_1}(\mu)} 
+ \Vert \nabla^j T\Vert_{L^{p_1}(\mu)}\right)\\
\,&\,\leq\,C\left(\Vert\nabla^{j+2} T\Vert_{L^{p_2}(\mu)} 
+ 2\,\Vert\nabla^{j+1} T\Vert_{L^{p_2}(\mu)} 
+ \Vert\nabla^{j} T\Vert_{L^{p_2}(\mu)}\right)\\
\,&\,\leq\qquad \dots\\
\,&\,\leq\,C\left(\Vert\nabla^{s} T\Vert_{L^{p_{s-j}}(\mu)} 
+ \dots + \Vert\nabla^j T\Vert_{L^{p_{s-j}}(\mu)}\right)\\
\,&\,\leq C\Vert T\Vert_{W^{s,p_{s-j}}(\mu)}\,.
\end{align*}
Since the $p_i$ are related by
$$
\frac{1}{p_{i}}=\frac{1}{p_{i+1}}-\frac{1}{n}\,,
$$
$p_0=p$ and $p_{s-j}=q$, we have
$$
\frac{1}{p}=\frac{1}{p_{s-j}}-\frac{s-j}{n}=\frac{1}{q}-\frac{s-j}{n}\,,
$$
and the first part of the corollary is proved.\\
The second part follows analogously using also
inequality~\eqref{sectens}.
\end{proof}

Now we put together this result and the universal inequalities
\begin{equation}\label{a=js}
    \Vert\nabla^j T\Vert_{L^{p}{(\mu)}}\leq\,C\,\Vert 
T\Vert_{W^{s,q}{(\mu)}}^{\frac{j}{s}}\Vert
      T\Vert_{L^r{(\mu)}}^{\frac{s-j}{s}}\,,
\end{equation}
which are obviously implied by Proposition~\ref{univ1gen}, to get the 
following interpolation type inequalities.

\begin{prop}
In the same hypothesis on $(M,g)$ as before, there exist a constant 
$C$ depending only on $n$, $l$, $s$, $j$, $p$, $q$, $r$, $\delta$ and
$B$, such that for every covariant tensor $T=T_{i_1 \dots i_l}$, the 
following inequality hold
\begin{equation}\label{gn}
    \Vert\nabla^j T\Vert_{L^{p}{(\mu)}}\leq\,C\,\Vert 
T\Vert_{W^{s,q}{(\mu)}}^{a}\Vert
      T\Vert_{L^r{(\mu)}}^{1-a}\,,
\end{equation}
for all $j\in[0,s]$, $p, q, r\in[1,+\infty)$ and $a\in[j/s,1]$
with the compatibility condition
$$
\frac{1}{p}=\frac{j}{n}+a\left(\frac{1}{q}-\frac{s}{n}\right)+\frac{1-a}{r}\,.
$$
If such condition gives a negative value for $p$, the inequality holds
for every $p\in[1,+\infty)$ on the left side.
\end{prop}

\begin{proof}
The cases $a=j/s$ and $a=1$ are inequalities~\eqref{a=js}
and~\eqref{a=1}, respectively, the intermediate cases, when $j/s<a<1$,
are obtained immediately by the $\log$--convexity of
$\Vert\,\cdot\,\Vert_{L^{p}(\mu)}$ in ${1/p}$, which is a linear
function of $a$, and the fact that the right side is exponential in
$a$.\\
If $p$ is negative then $\frac{1}{q}-\frac{s}{n}<0$ and
$$
\frac{1}{q}-\frac{s-j}{n}\leq
\frac{j}{n}+a\left(\frac{1}{q}-\frac{s}{n}\right)+\frac{1-a}{r}\,,
$$
hence, the $L^{\infty}$ estimate of inequality~\eqref{a=1inf} together
with~\eqref{a=js} gives the inequality for every $p\in[1,+\infty)$.
\end{proof}

\begin{rem}
By simplicity, we avoided to discuss in all the section the critical
cases of the inequalities, for instance $p=n$ in
Proposition~\ref{pnnp}. Actually, for our purposes, we just need to say
that in a critical case we can allow any value of $p\in[1,+\infty)$ 
in the left side of inequalities like~\eqref{gn}. This can be seen
easily, by considering a suitable inequality with a lower
integrability exponent on the right side and then applying H\"older
inequality.
\end{rem}

Putting together the estimates of this section with
Proposition~\ref{n+1} we obtain the following result.

\begin{prop} As long as the flow by the gradient of ${\mathcal F}_m$
of a hypersurface in $\Ri$ exists, for every smooth covariant tensor
$T=T_{i_1 \dots i_l}$ we have the inequalities
\begin{equation}\label{gnfinal}
    \Vert\nabla^j T\Vert_{L^{p}{(\mu)}}\leq\,C\,\Vert 
T\Vert_{W^{s,q}{(\mu)}}^{a}\Vert
      T\Vert_{L^r{(\mu)}}^{1-a}\,,
\end{equation}
for all $j\in[0,s]$, $p, q, r\in[1,+\infty)$ and $a\in[j/s,1]$
with the compatibility condition
$$
\frac{1}{p}=\frac{j}{n}+a\left(\frac{1}{q}-\frac{s}{n}\right)+\frac{1-a}{r}\,.
$$
If such condition gives a negative value for $p$, the inequality holds
for every $p\in[1,+\infty)$ on the left side.\\
The constant $C$ depends only on $m$, $n$, $l$, $s$, $j$, $p$, $q$,
$r$ and the value of ${\mathcal F}_m$ for the initial hypersurface.
\end{prop}

\section{Long Time Existence of the Flow}

Suppose that at a certain time $T>0$ the evolving hypersurface
develops a singularity, then considering the family
$\left\{M_t\right\}_{t\in[0,T)}$, we are going to use the time--independent
inequalities~\eqref{gnfinal} to show that we have uniform estimates
$$
\max_{M_t}\vert\nabla^k\AAA\vert\,\leq\,C_k\,<\,+\infty\qquad
\forall t\in[0,T)\,
$$
for all $k\in\NN$. We will see that such estimates are in
contradiction with the development of a singularity at time $t=T$,
hence the  flow must be smooth for every positive time.\\
To this aim we are going to study the evolution of the following
integrals,
$$
\int_{M}\vert\nabla^{k}\AAA\vert^2\,d\mu_t\,.
$$

\begin{rem} As in the previous sections, we will omit to say in the
  computations that all the polynomials $\pol_s$ and $\qol^s$
  which will appear are  independent of the manifold $(M,g)$ where the
  tensors are defined.
\end{rem}

First we derive the evolution equations for $g$, $\nu$,
$\Gamma_{jk}^i$ and $\AAA$. Essentially repeating the computations of
Section~\ref{firstvariat}, we get

\begin{align*}
\frac{\partial}{\partial t}g_{ij}\,
\,=&\,-2\EEE_mh_{ij}\\
\frac{\partial}{\partial t}g^{ij}\,
\,=&\, 2\EEE_mh^{ij}\\
\frac{\partial}{\partial t}\nu\phantom{^{ij}}\,
\,=&\,\phantom{-}\nabla \EEE_m\\
\frac{\partial}{\partial t}\Gamma_{jk}^i
\,=&\,\nabla\EEE_m * {\AAA} + \EEE_m * \nabla\AAA\,.
\end{align*}

\begin{lemma}
The second fundamental form of $M_t$ satisfies the 
evolution equation
$$
\frac{\partial}{\partial t}h_{ij}=
2(-1)^{m}\overset{\text{$m+1$~{\rm times}}}
{\overbrace{\Delta\comp\dots\comp\Delta}}h_{ij} + \qol^{2m+3}(\AAA,
\AAA) + \qol^{2m+3}(\nabla\nu, \AAA) + \qol^3(\AAA)\,.
$$
\end{lemma}

\begin{proof} Keeping in mind the Gauss--Weingarten
  relations~\eqref{gwein} and the equations above, we compute
\begin{align*}
\frac{\partial}{\partial t}h_{ij}
\,=&\,-\frac{\partial}{\partial t}\left\langle\nu\,\left\vert\,
\frac{\partial^2\varphi}{\partial x_i\partial
  x_j}\right.\right\rangle\\
\,=&\,\left\langle\nu\,\left\vert\,
\frac{\partial^2(\EEE_m\nu)}{\partial x_i\partial
  x_j}\right.\right\rangle-\left\langle\nabla\EEE_m\,\left\vert\,
\frac{\partial^2\varphi}{\partial x_i\partial  x_j}\right.\right\rangle\\
\,=&\,\frac{\partial^2\EEE_m}{\partial x_i\partial
  x_j}+\EEE_m\left\langle\nu\,\left\vert\,\frac{\partial}{\partial
    x_i}\left(h_{jl}g^{ls}\frac{\partial\varphi}{\partial
      x_s}\right)\right.\right\rangle\\
\,&\, - \left\langle\left.\frac{\partial\EEE_m}{\partial
    x_l}\cdot\frac{\partial\varphi}{\partial
    x_s}g^{ls}\,\right\vert\,\Gamma_{ij}^k\frac{\partial\varphi}{\partial x_k}
-h_{ij}\nu\right\rangle\\
\,=&\,\frac{\partial^2\EEE_m}{\partial x_i\partial
  x_j}-\Gamma_{ij}^k\frac{\partial \EEE_m}{\partial x_k}
+\EEE_mh_{jl}g^{ls}\left\langle\nu\,\left\vert\,
\Gamma_{is}^z\frac{\partial\varphi}{\partial
  x_z}-h_{is}\nu\right\rangle\right.\\
\,=&\,\nabla_i\nabla_j \EEE_m - \EEE_m h_{is}g^{sl}h_{lj}\,.
\end{align*}
Expanding $\EEE_m$ we continue,
\begin{align*}
\frac{\partial}{\partial t}h_{ij}
\,=&\,\nabla_i\nabla_j\Bigl(2(-1)^{m}\overset{\text{$m$~{\rm times}}}
{\overbrace{\Delta\Delta\dots\Delta}} \HHH + \qol^{2m+1}(\nabla\nu, \AAA) +
\qol^1(\AAA)\Bigr)\\
\,&\,-\Bigl(2(-1)^{m}\overset{\text{$m$~{\rm times}}}
{\overbrace{\Delta\Delta\dots\Delta}} \HHH + \qol^{2m+1}(\nabla\nu, \AAA) +
\qol^1(\AAA)\Bigr)h_{is}g^{sl}h_{lj}\\
\,=&\,2(-1)^m\,\nabla_i\nabla_j\overset{\text{$m$~{\rm times}}}
{\overbrace{\Delta\Delta\dots\Delta}} \HHH + \qol^{2m+3}(\nabla\nu,
\AAA) + \qol^3(\AAA)\,.
\end{align*}
Interchanging repeatedly derivatives in the first term we introduce
some extra terms of the form $\qol^{2m+3}(\AAA, \AAA)$ and we get
$$
\frac{\partial}{\partial t}h_{ij}
=2(-1)^m\,\overset{\text{$m$~{\rm times}}}
{\overbrace{\Delta\Delta\dots\Delta}}\nabla_i\nabla_j \HHH + 
\qol^{2m+3}(\AAA, \AAA) + \qol^{2m+3}(\nabla\nu, \AAA) +
\qol^3(\AAA)\,,
$$
then using equation~\eqref{codaz} we conclude
\begin{align*}
\frac{\partial}{\partial t}h_{ij}
\,=&\,2(-1)^m\,\overset{\text{$m$~{\rm times}}}
{\overbrace{\Delta\Delta\dots\Delta}}(\Delta h_{ij} - \HHH
h_{il}g^{ls}h_{sj} - \vert\AAA\vert^2h_{ij})\\
\,&\,+\qol^{2m+3}(\AAA, \AAA) + \qol^{2m+3}(\nabla\nu, \AAA) 
+ \qol^3(\AAA)\\
\,=&\,2(-1)^m\,\overset{\text{$m+1$~{\rm times}}}
{\overbrace{\Delta\Delta\dots\Delta}} h_{ij} +
\qol^{2m+3}(\AAA, \AAA) + \qol^{2m+3}(\nabla\nu, \AAA) +
\qol^3(\AAA)\,.
\end{align*}
\end{proof}

Now we deal with the covariant derivatives of $\AAA$.

\begin{lemma} We have
\begin{align*}
\frac{\partial}{\partial t}\nabla^k h_{ij}
\,=&\,2(-1)^m\overset{\text{$m+1$~{\rm times}}}
{\overbrace{\Delta\Delta\dots\Delta}} \nabla^kh_{ij}\\
\,&\,+ \qol^{k+2m+3}(\AAA, \AAA) + \qol^{k+2m+3}(\nabla\nu,
\AAA) + \qol^{k+3}(\AAA)\,.
\end{align*}
\end{lemma}

\begin{proof} With a reasoning analogous to the one of 
Lemma~\ref{eulerkey} applied to the tensor $\AAA$ and by the
  previous lemma, we have
\begin{align*}
\frac{\partial}{\partial t}\nabla^k h_{ij}
\,=&\,\nabla^k\frac{\partial}{\partial t}h_{ij} + \pol_k(\AAA,
\AAA, \EEE_m)\\
\,=&\,\nabla^k\frac{\partial}{\partial t}h_{ij} 
+ \qol^{k+2m+3}(\AAA, \AAA) + \qol^{k+2m+3}(\nabla\nu,
\AAA) + \qol^{k+3}(\AAA, \AAA)\\
\,=&\,2(-1)^m\nabla^k\overset{\text{$m+1$~{\rm times}}}
{\overbrace{\Delta\Delta\dots\Delta}} h_{ij}\\
\,&\,+\,\nabla^k\qol^{2m+3}(\AAA, \AAA) + \nabla^k\qol^{2m+3}(\nabla\nu,
\AAA) +\,\nabla^k\qol^3(\AAA)\\
\,&\,+ \qol^{k+2m+3}(\AAA, \AAA) + \qol^{k+2m+3}(\nabla\nu,
\AAA) + \qol^{k+3}(\AAA, \AAA)\\
\,=&\,2(-1)^m\nabla^k\overset{\text{$m+1$~{\rm times}}}
{\overbrace{\Delta\Delta\dots\Delta}} h_{ij}\\
\,&\,+ \qol^{k+2m+3}(\AAA, \AAA) + \qol^{k+2m+3}(\nabla\nu,
\AAA) + \qol^{k+3}(\AAA)\,.
\end{align*}
Interchanging the operator $\nabla^k$ with the Laplacians in the first
term and including the extra terms in $\qol^{k+2m+3}(\AAA, \AAA)$, we
obtain
\begin{align*}
\frac{\partial}{\partial t}\nabla^k h_{ij}
\,=&\,2(-1)^m\overset{\text{$m+1$~{\rm times}}}
{\overbrace{\Delta\Delta\dots\Delta}} \nabla^kh_{ij}\\
\,&\,+ \qol^{k+2m+3}(\AAA, \AAA) 
+ \qol^{k+2m+3}(\nabla\nu, \AAA) 
+ \qol^{k+3}(\AAA)\,.
\end{align*}
\end{proof}

\begin{prop} The following formula holds,
\begin{align*}
\frac{\partial}{\partial t}\int_M\vert\nabla^k\AAA\vert^2\,d\mu_t
\,=&\,-4\int_M \vert\nabla^{k+m+1}\AAA\vert^2\,d\mu_t\\
\,&\,+\int_M \qol^{2(k+m+2)}(\AAA, \AAA, \AAA) 
+ \qol^{2(k+m+2)}(\nabla\nu, \AAA, \AAA)\,d\mu_t\\ 
\,&\,+\int_M \qol^{2(k+2)}(\AAA, \AAA)\,d\mu_t\,.
\end{align*}
\end{prop}

\begin{proof}
By the previous results we have
\begin{align*}
\frac{\partial}{\partial t}\vert\nabla^k\AAA\vert^2
\,=&\,2g^{i_1j_1}\dots g^{i_kj_k}g^{is}g^{jz}\frac{\partial}{\partial
  t}\nabla_{i_1\dots i_k} h_{ij}\nabla_{j_1\dots j_k} h_{sz}\\
\,&\,+g^{i_1j_1}\dots\frac{\partial}{\partial t}g^{i_lj_l}\dots
g^{i_kj_k}g^{is}g^{jz}
\nabla_{i_1\dots i_k} h_{ij}\nabla_{j_1\dots j_k} h_{sz}\\
\,=&\,4(-1)^mg^{i_1j_1}\dots g^{i_kj_k}g^{is}g^{jz}
\overset{\text{$m+1$~{\rm times}}}
{\overbrace{\Delta\Delta\dots\Delta}}
\nabla_{i_1\dots i_k} h_{ij}\nabla_{j_1\dots j_k} h_{sz}\\
\,&\,+\Bigl(\qol^{k+2m+3}(\AAA, \AAA) + \qol^{k+2m+3}(\nabla\nu,
\AAA) + \qol^{k+3}(\AAA)\Bigr) * \nabla^k\AAA\\
\,&\,+2\EEE_mg^{i_1j_1}\dots h^{i_lj_l}\dots
g^{i_kj_k}g^{is}g^{jz}
\nabla_{i_1\dots i_k} h_{ij}\nabla_{j_1\dots j_k} h_{sz}\\
\,=&\,4(-1)^mg^{i_1j_1}\dots g^{i_kj_k}g^{is}g^{jz}
\overset{\text{$m+1$~{\rm times}}}
{\overbrace{\Delta\Delta\dots\Delta}}
\nabla_{i_1\dots i_k} h_{ij}\nabla_{j_1\dots j_k} h_{sz}\\
\,&\,+\qol^{2(k+m+2)}(\AAA, \AAA, \AAA) + \qol^{2(k+m+2)}(\nabla\nu,
\AAA, \AAA) + \qol^{2(k+2)}(\AAA, \AAA)\\
\,=&\,4(-1)^mg^{is}g^{jz} \nabla_{i_{k+1}}\nabla^{i_{k+1}}\dots
  \nabla_{i_{k+m+1}}\nabla^{i_{k+m+1}}\nabla_{i_1\dots i_k} h_{ij}
\nabla^{i_1\dots i_k} h_{sz}\\
\,&\,+\qol^{2(k+m+2)}(\AAA, \AAA, \AAA) + \qol^{2(k+m+2)}(\nabla\nu,
\AAA, \AAA) + \qol^{2(k+2)}(\AAA, \AAA)\,.
\end{align*}
Interchanging the covariant derivatives in the first term we introduce some
extra terms of the form $\qol^{2(k+m+2)}(\AAA, \AAA, \AAA)$, hence we get
\begin{align*}
\frac{\partial}{\partial t}&\,\int_M\vert\nabla^k\AAA\vert^2\,d\mu_t =\\
\,&\,4(-1)^m\int_M g^{is}g^{jz} \nabla^{i_{k+1}}\dots\nabla^{i_{k+m+1}}
\nabla_{i_{k+m+1}}\dots\nabla_{i_{k+1}}\nabla_{i_1\dots i_k} h_{ij}
\nabla^{i_1\dots i_k} h_{sz}\,d\mu_t\\
\,&\,+\int_M \qol^{2(k+m+2)}(\AAA, \AAA, \AAA) + \qol^{2(k+m+2)}(\nabla\nu,
\AAA, \AAA) + \qol^{2(k+2)}(\AAA, \AAA)\,d\mu_t\\
\,&\,+\int_M \qol^{2(k+2)}(\AAA, \AAA)\,d\mu_t\,,
\end{align*}
where the last integral comes from the time derivative of $\mu_t$.\\
Then, carrying the $m+1$ derivatives
$\nabla^{i_{k+1}}\dots\nabla^{i_{k+m+1}}$ on $\nabla^{i_1\dots i_k}
h_{sz}$ by means of the divergence theorem, we finally obtain the
claimed result,
\begin{align*}
\phantom{a}\,=&\,-4\int_M g^{is}g^{jz}
\nabla_{i_{k+m+1}}\dots\nabla_{i_{k+1}}\nabla_{i_1\dots i_k} h_{ij}
\nabla^{i_{k+m+1}}\dots\nabla^{i_{k+1}}\nabla^{i_1\dots i_k} h_{sz}\,d\mu_t\phantom{aaaaa}\\
\,&\,+\int_M \qol^{2(k+m+2)}(\AAA, \AAA, \AAA) + \qol^{2(k+m+2)}(\nabla\nu,
\AAA, \AAA) + \qol^{2(k+2)}(\AAA, \AAA)\,d\mu_t\\
\,=&\,-4\int_M \vert\nabla^{k+m+1}\AAA\vert^2\,d\mu_t\\
\,&\,+\int_M \qol^{2(k+m+2)}(\AAA, \AAA, \AAA) + \qol^{2(k+m+2)}(\nabla\nu,
\AAA, \AAA) + \qol^{2(k+2)}(\AAA, \AAA)\,d\mu_t\,.
\end{align*}
The leading coefficient became $-4$ since we
multiplied $4(-1)^m$ for $(-1)^{m+1}$ while doing the $m+1$
integrations by parts. 
\end{proof}

Now we analyze the terms
$$
\int_M \qol^{2(k+m+2)}(\AAA, \AAA, \AAA)\,d\mu_t\quad\text{ and }\quad
\int_M \qol^{2(k+m+2)}(\nabla\nu, \AAA, \AAA)\,d\mu_t\,.
$$
If one of the two polynomials contains a derivative $\nabla^i\AAA$ or 
$\nabla^{i}(\nabla\nu)$ of order $i>k+m+1$, then all the other
derivatives must be of order lower  than $k+m$, since the rescaling
order of the polynomials is $2(k+m+2)$ and the fact that there are at
least three factors in every additive term. 
In this case, using repeatedly the divergence theorem as before, to
lower such highest derivative, we get the integral of a new polynomial 
which does not contain derivatives of order higher than
$k+m+1$. Moreover, if there is a derivative of order $k+m+1$ then
the order of all the other derivatives in $\qol^{2(k+m+2)}$ must be
lower or equal than $k+m$, by the same argument.\\
With the same reasoning, the term
$$
\int_M \qol^{2(k+2)}(\AAA, \AAA)\,d\mu_t\,,
$$
can be transformed it in a term without derivatives of order higher or
equal than $k+m+1$.

Hence, we can suppose that the last three terms in
\begin{align}
\frac{\partial}{\partial t}\int_M\vert\nabla^k\AAA\vert^2\,d\mu_t
\,=&\,-4\int_M \vert\nabla^{k+m+1}\AAA\vert^2\,d\mu_t\nonumber\\
\,&\,+\int_M \qol^{2(k+m+2)}(\AAA, \AAA, \AAA) +
\qol^{2(k+m+2)}(\nabla\nu, \AAA, \AAA)\,d\mu_t\nonumber\\
\,&\,+\int_M \qol^{2(k+2)}(\AAA, \AAA)\,d\mu_t\label{stima999}
\end{align}
do not contain derivatives of $\AAA$ or of $\nabla\nu$ of order higher
than $k+m+1$; possibly, only one derivative of order $k+m+1$ can
appear.

\begin{lemma}\label{covacurv}
The following inequality holds
$$
\vert\nabla^s\nu\vert\leq
\vert\nabla^{s-1}\AAA\vert +  
\vert \qol^{s}(\AAA)\vert\,,
$$
where $\qol^{s}(\AAA)$ does not contain derivatives of $\AAA$ of order
higher than $s-2$.
\end{lemma}

\begin{proof} By equations~\eqref{gwein} it follows that
$\nabla\nu=\AAA * \nabla\varphi$, hence
$$
\nabla^s\nu=\nabla^{s-1}\AAA * \nabla\varphi + \sum_{i+j=s-2}
\nabla^i\AAA * \nabla^j\nabla^2\varphi
$$
and since $\nabla_{ij}^2\varphi=-h_{ij}\nu$, we get
\begin{align*}
\nabla^s\nu
\,=&\,\nabla^{s-1}\AAA * \nabla\varphi + \sum_{i+j=s-2}
\nabla^i\AAA * \nabla^j(\AAA\nu)\\
\,=&\,\nabla^{s-1}\AAA * \nabla\varphi + \sum_{i+j+k=s-2}
\nabla^i\AAA * \nabla^j\AAA * \nabla^k\nu\,.
\end{align*}
Then, by an induction argument we can express $\nabla^s\nu$ as 
$$
\nabla^s\nu=\nabla^{s-1}\AAA * \nabla\varphi + \qol^s(\AAA)
$$
where $\qol^s(\AAA)$ does not contain derivatives of order higher than 
$s-2$.\\
Taking the norm of both sides we get 
$$
\vert \nabla^s\nu\vert\leq
\vert\nabla^{s-1}\AAA * \nabla\varphi\vert 
+ \vert\qol^s(\AAA)\vert
$$
and we conclude the proof computing
\begin{align*}
\vert\nabla^{s-1}\AAA * \nabla\varphi\vert
\,=&\,\left\vert \nabla_{i_1\dots i_{s-1}}h_{il} 
g^{lk}\frac{\partial\varphi}{\partial x_k} \,\right\vert\\ 
\,=&\,\left(\nabla_{i_1\dots i_{s-1}}h_{il} 
g^{lk}\frac{\partial\varphi}{\partial x_k} 
g^{i_1j_1}\dots g^{i_{s-1}j_{s-1}}g^{ij}
\nabla_{j_1\dots j_{s-1}}h_{jw} 
g^{wz}\frac{\partial\varphi}{\partial x_z}\right)^{1/2}\\
\,=&\,\left(\nabla_{i_1\dots i_{s-1}}h_{il} 
g^{lk}g_{kz}g^{wz}
g^{i_1j_1}\dots g^{i_{s-1}j_{s-1}}g^{ij}
\nabla_{j_1\dots j_{s-1}}h_{jw}\right)^{1/2}\\
\,=&\,\left(\nabla_{i_1\dots i_{s-1}}h_{il} 
g^{lw}g^{i_1j_1}\dots g^{i_{s-1}j_{s-1}}g^{ij}
\nabla_{j_1\dots j_{s-1}}h_{jw}\right)^{1/2}\\
\,=&\,\vert\nabla^{s-1}\AAA\vert\,. 
\end{align*}
\end{proof}

Taking the absolute values inside the integrals and using this lemma
to substitute every derivative of $\nu$ in~\eqref{stima999}, we obtain
$$
\frac{\partial}{\partial t}\int_M\vert\nabla^k\AAA\vert^2\,d\mu_t
\leq-4\int_M \vert\nabla^{k+m+1}\AAA\vert^2\,d\mu_t + 
\int_M \vert\qol^{2(k+m+2)}(\AAA)\vert + 
\vert\qol^{2(k+2)}(\AAA)\vert\,d\mu_t
$$
where, as before, the two polynomials do not contain derivatives of
$\AAA$ of order higher than $k+m+1$; possibly, only one derivative of
order $k+m+1$ can appear in every multiplicative term of
$\qol^{2(k+m+2)}(\AAA)$.

Before going on, we remark that the $*$ product of
tensors satisfies the following metric property,
\begin{equation}\label{normstima}
\vert T * S\vert\leq \vert T\vert \cdot \vert S\vert\,.
\end{equation}
This can be easily seen choosing an orthonormal basis at a point of
$M$, in such coordinates we have
\begin{align*}
\vert T * S\vert^2
\,=&\,\sum_{\genfrac{}{}{0pt}{}{\text{free}}{\text{indices}}}
\biggl(\sum_{\genfrac{}{}{0pt}{}{\text{contracted}}{\text{indices}}}
T_{i_1\dots i_k}S_{j_1\dots j_l}\biggr)^2\\
\,\leq&\,\sum_{\genfrac{}{}{0pt}{}{\text{free}}{\text{indices}}}
\biggl(\sum_{\genfrac{}{}{0pt}{}{\text{contracted}}{\text{indices}}} T_{i_1\dots i_k}^2\biggr) 
\biggl(\sum_{\genfrac{}{}{0pt}{}{\text{contracted}}{\text{indices}}} S_{j_1\dots j_l}^2\biggr)\\
\,\leq&\,
\biggl(\sum_{\genfrac{}{}{0pt}{}{\text{free}}{\text{indices}}}
\sum_{\genfrac{}{}{0pt}{}{\text{contracted}}{\text{indices}}} 
T_{i_1\dots i_k}^2\biggr) 
\biggl(\sum_{\genfrac{}{}{0pt}{}{\text{free}}{\text{indices}}}
\sum_{\genfrac{}{}{0pt}{}{\text{contracted}}{\text{indices}}} 
S_{j_1\dots j_l}^2\biggr)\\
\,=&\,\vert T \vert^2 \cdot \vert S\vert^2\,.
\end{align*}

Now by definition we have
$$
\qol^{2(k+m+2)}(\AAA)=\sum_j \bbigstar_{l=1}^{N_j}\nabla^{c_{jl}}\AAA
$$
with
$$
\sum_{l=1}^{N_j} (c_{jl}+1)=2(k+m+2)
$$
for every $j$, hence
$$
\vert\qol^{2(k+m+2)}(\AAA)\vert\leq\sum_j \prod_{l=1}^{N_j}\vert\nabla^{c_{jl}}\AAA\vert
$$
by~\eqref{normstima}.
Setting
$$
Q_j=\prod_{l=1}^{N_j}\vert\nabla^{c_{jl}}\AAA\vert
$$
we clearly obtain
$$
\int_M \vert\qol^{2(k+m+2)}(\AAA)\vert\,d\mu_t\leq\sum_j\int_M Q_j\,d\mu_t\,.
$$
If $Q_j$ contains a derivative of $\AAA$ of order $k+m+1$, we have
seen that all the others have order lower or equal than $k+m$, then
collecting derivatives of the same order,  $Q_j$ can be estimated as
follows
$$
Q_j\leq\vert\nabla^{k+m+1}\AAA\vert\, \cdot
\prod_{i=0}^{k+m}\vert\nabla^i\AAA\vert^{\alpha_{ji}}
$$
for some $\alpha_{ji}$ satisfying the rescaling condition
$$
(k+m+2)+\sum_{i=0}^{k+m}(i+1)\alpha_{ji}=2(k+m+2)\,.
$$
Hence, using Young inequality, for every $\varepsilon_j>0$ we have
\begin{align*}
\int_M Q_j\,d\mu_t
\,\leq&\,\varepsilon_j\int_M \vert\nabla^{k+m+1}\AAA\vert^2\,d\mu_t
+ \frac{1}{4\varepsilon_j}\int_M
\prod_{i=0}^{k+m}\vert\nabla^i\AAA\vert^{2\alpha_{ji}}\,d\mu_t\\
\,=&\,\varepsilon_j\int_M \vert\nabla^{k+m+1}\AAA\vert^2\,d\mu_t
+ \int_M \vert\qol^{2(k+m+2)}(\AAA)\vert\,d\mu_t\,,
\end{align*}
where we put in evidence the fact that the last term satisfies again
the rescaling condition and no more contains the derivative
$\nabla^{k+m+1}\AAA$.\\
Collecting all together such ``bad'' terms, and choosing suitable
$\varepsilon_j>0$ such that their total sum is less than one, 
we obtain
$$
\frac{\partial}{\partial t}\int_M\vert\nabla^k\AAA\vert^2\,d\mu_t
\leq-3\int_M \vert\nabla^{k+m+1}\AAA\vert^2\,d\mu_t + 
\int_M \vert\qol^{2(k+m+2)}(\AAA)\vert + 
\int_M \vert\qol^{2(k+2)}(\AAA)\vert\,d\mu_t
$$
where now in the last two terms all the derivatives of $\AAA$ have
order lower than $k+m+1$. We are then ready to estimate them via
interpolation inequalities.

As before,  
$$
\vert\qol^{2(k+m+2)}(\AAA)\vert\leq\sum_j Q_j
$$
and after collecting derivatives of the same order in $Q_j$,
$$
Q_j = \prod_{i=0}^{k+m}\vert\nabla^i \AAA\vert^{\alpha_{ji}} \qquad
\text{ with } \quad \sum_{i+1}^{k+m}\alpha_{ji}(i+1) = 2(k+m+2)\,.
$$
Then,
\begin{align*}
\int_{M} Q_j\,d\mu_t=\,&\,
\int_{M}\prod_{i=0}^{k+m}\vert\nabla^i
\AAA\vert^{\alpha_{ji}}\,d\mu_t\\
\leq\,&\,\prod_{i=0}^{k+m}\left(\int_{M}\vert\nabla^i
    \AAA\vert^{\alpha_{ji}\gamma_{i}}\,d\mu_t\right)^{\frac{1}{\gamma_i}}\\
=\,&\,\prod_{i=0}^{k+m}\Vert\nabla^i
\AAA\Vert_{L^{\alpha_{ji}\gamma_i}(\mu_t)}^{\alpha_{ij}}
\end{align*}
where the $\gamma_i$ are arbitrary positive values such that $\sum
1/{\gamma_i}=1$.

We apply interpolation inequalities: if in~\eqref{gn} we take $q=2$,
$r=n+1$, $s=k+m+1$, $j=i$ and $T=\AAA$ we get
$$
    {\Vert\nabla^{i} \AAA\Vert}_{L^{p_i}(\mu_t)}
\leq C{\Vert\AAA\Vert}_{W^{2,k+m+1}(\mu_t)}^{a}
{\Vert\AAA\Vert}_{L^{n+1}(\mu_t)}^{1-a}
$$
with
\begin{equation}\label{ccc9}
a=\frac{\frac{1}{p_i}-\frac{i}{n}-\frac{1}{n+1}}
{\frac{1}{2}-\frac{k+m+1}{n}-\frac{1}{n+1}}\in\left[\frac{i}{k+m+1},1\right]
\end{equation}
and $p_i>1$.\\
Now, since the volumes of $M_t$ and
${\Vert\AAA\Vert}_{L^{n+1}(\mu_t)}$ are uniformly bounded in time, also
${\Vert\AAA\Vert}_{L^{2}(\mu_t)}$ is uniformly bounded and using the
universal inequalities~\eqref{a=js} with $p=q=r=2$ we have
\begin{align*}
{\Vert\AAA\Vert}_{W^{2,k+m+1}(\mu_t)}
\,\leq&\,\sum_{s=0}^{k+m+1}
C{\Vert\nabla^{k+m+1}\AAA\Vert}_{L^2(\mu_t)}^{\frac{s}{k+m+1}}\\
\,\leq&\,\sum_{s=0}^{k+m+1}
{\Vert\nabla^{k+m+1}\AAA\Vert}_{L^2(\mu_t)}+C\\
\,\leq&\,B{\Vert\nabla^{k+m+1}\AAA\Vert}_{L^2(\mu_t)}+C\,,
\end{align*}
where we applied Young inequality.\\
Hence, we conclude that we have constants $B$, $C$ independent of  
$t$ such that
\begin{equation}\label{gn20}
    {\Vert\nabla^{i} \AAA\Vert}_{L^{p_i}(\mu_t)}\leq \left(B{\Vert \nabla^{k+m+1}
      \AAA\Vert}_{L^2(\mu_t)}+C\right)^{a}
\end{equation}
for $a$ as in~\eqref{ccc9} and $p_i>1$.

Choosing $\gamma_i=0$ if $\alpha_{ji}=0$ and
$\gamma_i=\frac{2(k+m+2)}{\alpha_{ji}(i+1)}$ otherwise, we have clearly
$$
\sum_{i=0}^{k+m}\frac{1}{\gamma_i}=
\sum_{i=0}^{k+m}\frac{\alpha_{ji}(i+1)}{2(k+m+2)}=1\,
$$
by the rescaling condition on the $\alpha_{ji}$.\\
We claim that for every $i\in\{0, \dots , k+m\}$, the product
$p_i=\alpha_{ji}\gamma_i$ satisfies the condition~\eqref{ccc9}.\\
By definition, $p_i=\frac{2(k+m+2)}{i+1}$, hence we must check that
the following inequality holds
\begin{equation*}
\frac{i}{k+m+1}\leq\frac{\frac{i+1}{2(k+m+2)}-\frac{i}{n}-\frac{1}{n+1}}
{\frac{1}{2}-\frac{k+m+1}{n}-\frac{1}{n+1}}\leq 1
\end{equation*}
for every $i\in\{0, \dots , k+m\}$. Since every term is an affine
function of $i$, the claim follows if we show that the inequality holds
for $i=0$ and  $i=k+m+1$.\\
If $i=0$ we have to prove that
\begin{equation*}
0\leq\frac{\frac{1}{2(k+m+2)}-\frac{1}{n+1}}
{\frac{1}{2}-\frac{k+m+1}{n}-\frac{1}{n+1}}\leq 1\,,
\end{equation*}
that is, since the denominator of the fraction is negative 
(as $2m\geq n+1$),
\begin{equation*}
\frac{1}{2}-\frac{k+m+1}{n}-\frac{1}{n+1}\leq\frac{1}{2(k+m+2)}-\frac{1}{n+1}
\leq 0\,.
\end{equation*}
The right inequality is clearly true, again since $2m\geq n+1$, the
left one becomes
\begin{equation*}
\frac{k+m+1}{2(k+m+2)}=\frac{1}{2}-\frac{1}{2(k+m+2)}\leq\frac{k+m+1}{n}
\end{equation*}
which is true as $2(k+m+2)\geq n$.\\
When $i=k+m+1$ the fraction is equal to 1, hence the
inequality obviously holds. 

Then, the exponents $p_i=\alpha_{ji}\gamma_i$ are allowed in
inequality~\eqref{gn20} and we get
$$
\Vert\nabla^i\AAA\Vert_{L^{\alpha_{ji}\gamma_i}(\mu_t)}\leq
\left(B\Vert\nabla^{k+m+1}\AAA\Vert_{L^2(\mu_t)}+C\right)^{a_{ji}}
$$
where $a_{ji}$ is the relative value we obtain
from~\eqref{ccc9}.

Hence, 
\begin{align*}
\int_{M} Q_j\,d\mu_t
\leq\,&\,\prod_{i=0}^{k+m}\Vert\nabla^i
\AAA\Vert_{L^{\alpha_{ji}\gamma_i}(\mu_t)}^{\alpha_{ij}}\\
\leq\,&\,\prod_{i=0}^{k+m}\left(B\Vert\nabla^{k+m+1}
  \AAA\Vert_{L^2(\mu_t)}+ C\right)^{ a_{ji}\alpha_{ji}}\\
\leq\,&\,\left(B\Vert\nabla^{k+m+1}
  \AAA\Vert_{L^2(\mu_t)}+C\right)^{\sum_{i=0}^{k+m} a_{ji}\alpha_{ji}}
\end{align*}
where the constants $B$ and $C$ are independent of $t$ and 
$$
a_{ji}=\frac{\frac{1}{\alpha_{ji}\gamma_i}-\frac{i}{n}-\frac{1}{n+1}}{\frac{1}{2}-\frac{k+m+1}{n}-\frac{1}{n+1}}\,.
$$
Multiplying this relation by $\alpha_{ji}$ and summing on $i$ from $0$ to
$k+m$ we get
\begin{align*}
\sum_{i=0}^{k+m}
\alpha_{ji} a_{ji}\,=&\,\sum_{i=0}^{k+m}\frac{\frac{1}{\gamma_i}-\frac{i\alpha_{ji}}
{n}-\frac{\alpha_{ji}}{n+1}}{\frac{1}{2}-\frac{k+m+1}{n}-\frac{1}{n+1}}\\
=&\,\frac{1-\sum_{i=0}^{k+m}
\left(\frac{i\alpha_{ji}}{n}+\frac{\alpha_{ji}}{n+1}\right)}
{\frac{1}{2}-\frac{k+m+1}{n}-\frac{1}{n+1}}\\
=&\,\frac{1-\sum_{i=0}^{k+m}\frac{\alpha_{ji}(i+1)}{n}
-\sum_{i=0}^{k+m}{\alpha_{ji}}\left(\frac{1}{n+1}-\frac{1}{n}\right)}
{\frac{1}{2}-\frac{k+m+1}{n}-\frac{1}{n+1}}
\end{align*}
recalling that $\sum_{i=0}^{k+m}\alpha_{ji}(i+1)=2(k+m+2)$ we
continue,
\begin{align*}
\phantom{\sum_{i=0}^{k+m}\alpha_{ji} a_{ji}}
=&\,\frac{1-2\frac{k+m+2}{n}+\sum_{i=0}^{k+m}\frac{\alpha_{ji}}{n(n+1)}}
{\frac{1}{2}-\frac{k+m+1}{n}-\frac{1}{n+1}}\\
=&\,\frac{1-2\frac{k+m+1}{n}-\frac{2}{n}+\sum_{i=0}^{k+m}\frac{\alpha_{ji}}{n(n+1)}}
{\frac{1}{2}-\frac{k+m+1}{n}-\frac{1}{n+1}}\,.
\end{align*}
Now, the denominator is negative and clearly
$$
\sum_{i=0}^{k+m}{\alpha_{ji}}\geq\sum_{i=0}^{k+m}\frac{\alpha_{ji}(i+1)}{k+m+1}
=2\frac{k+m+2}{k+m+1}\,, 
$$
so we obtain
\begin{align*}
\sum_{i=0}^{k+m}
\alpha_{ji}
a_{ji}\,\leq
&\,\frac{1-2\frac{k+m+1}{n}-\frac{2}{n}+2\frac{k+m+2}{k+m+1}\frac{1}{n(n+1)}}
{\frac{1}{2}-\frac{k+m+1}{n}-\frac{1}{n+1}}\\
=&\,\frac{1-2\frac{k+m+1}{n}-\frac{2}{n}+\frac{2}{n(n+1)}
+\frac{2}{k+m+1}\frac{1}{n(n+1)}}
{\frac{1}{2}-\frac{k+m+1}{n}-\frac{1}{n+1}}\\
=&\,\frac{1-2\frac{k+m+1}{n}-\frac{2}{n+1}+\frac{2}{k+m+1}\frac{1}{n(n+1)}}
{\frac{1}{2}-\frac{k+m+1}{n}-\frac{1}{n+1}}\\
=&\,2-\frac{\frac{2}{k+m+1}\frac{1}{n(n+1)}}
{\frac{k+m+1}{n}+\frac{1}{n+1}-\frac{1}{2}}\\
=&\,2-\frac{4}{(k+m+1)[2(k+m+1)(n+1)-n(n-1)]}\, < \,2\,.
\end{align*}

Hence, we finally get
$$
\int_{M} Q_j\,d\mu_t\leq\left(B\int_{M}\vert\nabla^{k+m+1}
\AAA\vert^2\,d\mu_t+C\right)^{1-\delta}
$$
for a positive $\delta$ and using again Young inequality, we have
$$
\int_{M} Q_j\,d\mu_t\leq\varepsilon_j\int_{M}\vert\nabla^{k+m+1}
\AAA\vert^2\,d\mu_t+C
$$
for arbitrarily small $\varepsilon_j$. Repeating this argument for all
the $Q_j$ and choosing suitable $\varepsilon_j$ whose sum is less than
one, we conclude that
$$
\frac{d}{dt}\int_{M}\vert\nabla^k\AAA\vert^2\,\mu_t\,\leq
\,-2\int_{M}\vert\nabla^{k+m+1}\AAA\vert^2\,\mu_t+C+
\int_M\vert\qol^{2(k+2)}(\AAA)\vert\,d\mu_t
$$
with a constant $C$ independent of time.

The last term can be treated in the same way. 
It can be estimated by the sum of the multiplicative terms $Q_j$ 
and collecting derivatives of the same order as before, we have
$$
Q_j\leq\prod_{i=0}^{k+m}\vert\nabla^i \AAA\vert^{\beta_{ji}}\quad \text { with } \quad
\sum_{i=0}^{k+m}\beta_{ji}(i+1) = 2k+4\,.
$$
In this case the coefficients $\gamma_i$, when $\beta_{ji}\not=0$, are
given by  $\gamma_i=\frac{2(k+2)}{\alpha_{ji}(i+1)}$, hence
$$
\sum_{i=0}^{k+m}\frac{1}{\gamma_i}=\sum_{i=0}^{k+m}
\frac{\alpha_{ji}(i+1)}{2(k+2)}=1\,
$$
by the rescaling condition.\\
With an analogous control, one can see that the
conditions on the exponent $p_i$ are satisfied. It lasts to compute
\begin{align*}
\sum_{i=0}^{k+m}
\beta_{ji} a_{ji}\,=&\,\sum_{i=0}^{k+m}\frac{\frac{1}{\gamma_i}-\frac{i\beta_{ji}}
{n}-\frac{\beta_{ji}}{n+1}}{\frac{1}{2}-\frac{k+m+1}{n}-\frac{1}{n+1}}\\
=&\,\frac{1-\sum_{i=0}^{    k+m}\left(\frac{i\beta_{ji}}{n}+\frac{\beta_{ji}}{n+1}\right)}
{\frac{1}{2}-\frac{k+m+1}{n}-\frac{1}{n+1}}\\
=&\,\frac{1-\sum_{i=0}^{    k+m}\frac{\beta_{ji}(i+1)}{n}+\sum_{i=0}^{
k+m}\frac{\beta_{ji}}{n(n+1)}}
{\frac{1}{2}-\frac{k+m+1}{n}-\frac{1}{n+1}}\\
=&\,\frac{1-\frac{2k+4}{n}+\sum_{i=0}^{
k+m}\frac{\beta_{ji}}{n(n+1)}}
{\frac{1}{2}-\frac{k+m+1}{n}-\frac{1}{n+1}}\,.
\end{align*}
As the denominator is negative and
$$
\sum_{i=0}^{k+m}{\beta_{ji}}\geq\sum_{i=0}^{
k+m}\frac{\beta_{ji}(i+1)}{k+m+1}=\frac{2k+4}{k+m+1}\,,
$$
we obtain
\begin{align*}
\sum_{i=0}^{k+m}
\beta_{ji} a_{ji}\,\leq&\,\frac{1-\frac{2k+4}{n}+\sum_{i=0}^{
    k+m}\frac{\beta_{ji}(i+1)}{k+m+1}\frac{1}{n(n+1)}}
{\frac{1}{2}-\frac{k+m+1}{n}-\frac{1}{n+1}}\\
=&\,\frac{1-\frac{2k+4}{n}+\frac{2k+4}{k+m+1}\frac{1}{n(n+1)}}
{\frac{1}{2}-\frac{k+m+1}{n}-\frac{1}{n+1}}\,<2\,,
\end{align*}
since this last inequality is equivalent to
$$
1-\frac{2k+4}{n}+\frac{2k+4}{k+m+1}\frac{1}{n(n+1)}>
1-\frac{2(k+m+1)}{n}-\frac{2}{n+1}
$$
and simplifying, to 
$$
\frac{2k+4}{k+m+1}\frac{1}{n(n+1)}>
-\frac{2(m-1)}{n}-\frac{2}{n+1}
$$
which is obviously true.

Concluding as before we finally get
\begin{equation}\label{ddd10}
\frac{d}{dt}\int_{M}\vert\nabla^k\AAA\vert^2\,\mu_t\,\leq
\,-\int_{M}\vert\nabla^{k+m+1}\AAA\vert^2\,\mu_t+C
\end{equation}
for a constant $C$ independent of time.\\
By~\eqref{ccc3} and Young inequality, we have
\begin{align*}
\int_{M}\vert\nabla^k \AAA\vert^2\,\mu_t+C\,
&\leq\,B\Vert\nabla^{k+m+1}\AAA\Vert^{\frac{k}{k+m+1}}_{L^2{(\mu_t)}}
\Vert\AAA\Vert_{L^2{(\mu_t)}}^{\frac{m+1}{k+m+1}} + C\\
&\leq\,B\Vert\nabla^{k+m+1}\AAA\Vert^{\frac{k}{k+m+1}}_{L^2{(\mu_t)}}
+ C\\
&\leq\,\frac{1}{2}\int_{M}\vert\nabla^{k+m+1}
\AAA\vert^2\,\mu_t+C
\end{align*}
again with a uniform constant. Combining this inequality
with~\eqref{ddd10}, we obtain
$$
\frac{d}{dt}\int_{M}\vert\nabla^k\AAA\vert^2\,\mu_t\,\leq
\,-\,\frac{1}{2} \int_{M}\vert\nabla^{k}\AAA\vert^2\,\mu_t + C
$$
and a simple ODE's argument proves that there exists constants $C_k$ 
independent of time such that
$$
\int_{M}\vert\nabla^k \AAA\vert^2\,d\mu_t\leq C_k\,.
$$

To pass from $W^{2,p}(\mu_t)$ to pointwise estimates, first we
notice that being all the derivatives of $\AAA$ bounded 
in $L^2(\mu_t)$, by inequalities~\eqref{firstens}, for every $p\geq1$
and $k\in\NN$ we have constants $C_{k,p}$ such that
$$
\int_{M}\vert\nabla^k \AAA\vert^p\,d\mu_t\leq C_{k,p}\,.
$$
Then choosing a $p>n$, we apply inequalities~\eqref{sectens} to every
$\nabla^k\AAA$ to conclude that for every $k\in\NN$ we have constants $C_k$,
independent of $t$, such that
\begin{equation}\label{bbbbb}
\max_{M_t} \vert \nabla^k\AAA\vert \leq C_k\,.
\end{equation}

Looking back at the way we obtained them, 
we can see that the constants $C_k$ depend only on the
dimension $n$, the differentiation order $k$ and the initial
hypersurface $\varphi_0$.

Following Huisken~\cite{huisk1}, Sec.~8 and Kuwert and
Sch\"atzle~\cite{kuschat1}, Sec.~4, these estimates imply the
smoothness of the map $\varphi(p,t)$.\\
Since $\nabla^k\AAA$ are uniformly bounded in time, supposing that
$[0,T)$ is the maximal interval of existence of the flow, we have
$$
\vert\varphi(p,t)-\varphi(p,s)\vert\leq\int_s^t
\vert\EEE_m(\varphi_\xi)(p)\vert\,d\xi\leq C(t-s)
$$
for every $0\leq s\leq t<T$, then $\varphi_t$ uniformly converge to a 
continuous limit $\varphi_T$ as $t\to T$.\\
We recall Lemma~8.2 in~\cite{huisk1} (Lemma~14.2 in~\cite{hamilton1}).
\begin{lemma}
Let $g_{ij}$ a time--dependent metric on a compact manifold $M$ for
$0\leq t < T\leq +\infty$. Suppose that
$$
\int_0^T\max_{M_t}\left\vert\frac{\partial}{\partial
    t}g_{ij}\right\vert\,dt\leq C\,.
$$
Then the metrics $g_{ij}(t)$ are all equivalent, and they converge as
$t\to T$ uniformly to a positive definite metric tensor $g_{ij}(T)$
  which is continuous and also equivalent.
\end{lemma}
 
In our situation, if $T<+\infty$, the hypotheses of this lemma are
clearly satisfied, hence $\varphi(\cdot, T)$ represents a hypersurface.
Moreover, it also follows that there exists a positive constant $C$ 
depending only on $n$ and $\varphi_0$ 
such that for every $0\leq t<T$ we have
\begin{equation*}
\frac{1}{C}\leq g_{ij}(t)\leq C\,.
\end{equation*}
Since
$$
\frac{\partial}{\partial t}g_{ij} = -2\EEE_m h_{ij}
$$
by~\eqref{bbbbb}, for every $k\in\NN$ we have
$$
\left\Vert\nabla^k\frac{\partial}{\partial
t}g_{ij}\right\Vert_{L^\infty(\mu)} \leq C_k\,,
$$
analogously, as the time derivative of the Christoffel symbols is given 
by
$$
\frac{\partial}{\partial t}\Gamma_{jk}^i = \nabla\EEE_m * {\AAA} +
\EEE_m * \nabla\AAA
$$
it follows that
$$
\left\Vert\nabla^k\frac{\partial}{\partial
t}\Gamma_{jk}^i\right\Vert_{L^\infty(\mu)} \leq C_k\,.
$$
for every $k\in\NN$.\\
With an induction argument, we can prove the following formula
(where we avoid to indicate the indices) relating the iterated
covariant and coordinate derivatives of a tensor $T$,
\begin{equation}\label{coordin}
\nabla^m T=\partial^m T+\sum_{i=1}^m\sum_{j_1 + \dots + j_i+k\leq m-1} 
\partial^{j_1}\Gamma \dots \partial^{j_i}\Gamma \partial^k T\,.
\end{equation}
By this formula and induction, it follows that
$$
\Vert\partial^k\Gamma_{jl}^i\Vert_{L^\infty(\mu)}\,,\qquad
\left\Vert\partial^k\frac{\partial}{\partial
  t}\Gamma_{jl}^i\right\Vert_{L^\infty(\mu)} \leq C_k\,,
$$
for every $t\in[0,T)$.\\
Applying again formula~\eqref{coordin} to $T=\nabla^s\AAA$ we see that
$$
\partial^k\nabla^s\AAA-\nabla^{k+s}\AAA=
\sum_{i=1}^k\sum_{j_1 + \dots + j_i+l\leq k-1} 
\partial^{j_1}\Gamma \dots \partial^{j_i}\Gamma \, \partial^l
\nabla^s\AAA\,,
$$
and by induction and estimates~\eqref{bbbbb} we obtain
\begin{equation*}
\Vert\partial^k\nabla^s\AAA\Vert_{L^\infty(\mu)} \leq C_{k,s}
\end{equation*}
for every $k, s\in\NN$.\\
Since we already know that $\vert\varphi\vert$ is bounded and 
$\vert \partial\varphi\vert=1$, by the Gauss--Weingarten 
relations~\eqref{gwein}
$$
{\partial^2\varphi}=\Gamma\partial\varphi 
+ \AAA \nu\,,\qquad \partial\nu=\AAA * \partial\varphi
$$
and the previous estimates, we can conclude that
$$
\Vert\partial^k\varphi\Vert_{L^\infty(\mu)} \leq C_k
$$
for every $k\in\NN$ and $t\in[0,T)$.\\
The regularity of the time derivatives also follows by these 
estimates and the evolution equation.

Hence, the convergence $\varphi_t\to\varphi_T$, 
when $t\to T$, is in the $C^{\infty}$ topology and $M_T$ is
smooth. Then, using Theorem~\ref{smalltime} to restart the flow with
$\varphi_T$ as initial hypersurface, we get a contradiction with the
fact that $[0,T)$ is the maximal interval of  existence.

\begin{rem}
Though this argument shows that the solution is
classical, we cannot conclude that the estimates on the
parametrization hold uniformly for every $t\in[0,+\infty)$ which is
instead the case for the estimates~\eqref{bbbbb} on the curvature.
\end{rem}

\begin{teo}
If $m>\left[\frac{n}{2}\right]$, for any smooth hypersurface immersion
$\varphi_0:M\to\Ri$ there exists a unique smooth solution to the
problem
$$
\frac{\partial\varphi}{\partial t}(p,t)=-\EEE_m(\varphi_t)(p)\nu(p,t)\,,
$$
that is, the gradient flow associated to the functional 
$$
{\mathcal F}_m(\varphi)=\int_M1+\vert\nabla^m\nu\vert^2\,d\mu\,,
$$
defined for every $t\in[0,+\infty)$ and taking $\varphi_0$ as its
initial value.\\
Moreover, such solution satisfies
$$ 
\max_{M_t} \vert \nabla^k\AAA\vert \leq C_k\,.
$$
for constants $C_k$ depending only on $n$, $k$ and $\varphi_0$.
\end{teo}

\section{Convergence}

Let us consider the function $\sigma:[0,+\infty)\to\R$,
$$
\sigma(t)=\int_M\left[\EEE_m(\varphi_t)\right]^2\,d\mu_t\geq0\,.
$$
Clearly we have
$$
\frac{d}{dt}{\mathcal
  F}_m(\varphi_t)=-\int_M\left[\EEE_m(\varphi_t)\right]^2\,d\mu_t=-\sigma(t)\,,
$$
and integrating both sides in $t$ on $[0,+\infty)$ we get
$$
\int_0^{+\infty}\sigma(t)\,dt={\mathcal F}_m(\varphi_0)-{\mathcal
  F}_m(\varphi_t)\leq {\mathcal F}_m(\varphi_0)\,.
$$
Moreover,
$$
\left\vert\frac{d}{dt}\sigma(t)\right\vert=
\int_M\left\vert 2 \, \frac{\partial\EEE_m(\varphi_t)}
{\partial t} \, \EEE_m(\varphi_t) - \HHH 
\left[\EEE_m(\varphi_t)\right]^3\,\right\vert\,d\mu_t\leq C
$$
by the bounds~\eqref{bbbbb}. Then the function $\sigma$, being
Lipschitz and integrable on $[0,+\infty)$, converges
to zero at $+\infty$. This means that every $C^\infty$ limit 
hypersurface of the flow $\psi:M\to\Ri$ satisfies $\EEE_m(\psi)=0$,
i.\,e., it is a critical point of ${\mathcal F}_m$.

To find limit hypersurfaces, we need the following compactness
result of Langer and Delladio~\cite{della1,langer2}.

\begin{teo}\label{landel}
Let be given a family $(M,g_i)$ of closed, oriented,
$n$--dimensional hypersurfaces, isometrically immersed in $\Ri$ via
the maps $\varphi_i:M\to\Ri$, let $\mu_i$ the
associated measures on $M$ and
${\mathrm {Bar}}_i$ the {\em center of gravity} of $\varphi_i$,
that is,
$$
\text{${\mathrm {Bar}}$}_i=\int_M\varphi_i\,d\mu_i\,.
$$
Let $h$ be any metric tensor on $M$, 
if for some exponent $p>n$ and $C>0$ we have
$$
\int_M1+\vert\AAA\vert^p\,d\mu_i+\vert\text{${\mathrm
    {Bar}}$}_i\vert\,\leq\,C<+\infty\,,
$$
then there exist a subsequence of $\{\varphi_i\}$ (not
relabeled) and diffeomorphisms $\sigma_i: M\to M$ such that, 
$\{\varphi_i\circ\sigma_i\}$ converges in the $H^{2,p}$ weak topology
of maps from $(M,h)\to\Ri$ to an immersion $\varphi:M\to\Ri$.
\end{teo}

Translating the hypersurfaces $\varphi_t:M\to\R$ in order
to have Bar$_t=0\in\Ri$, we are in the above hypotheses. Hence, we can
extract a subsequence of smooth hypersurfaces $\varphi_i=\varphi_{t_i}$
and diffeomorphisms $\sigma_i: M\to M$ such that, for a fixed metric $h$
on $M$, the sequence  $\{\varphi_i\circ\sigma_i\}$ converges in the
$H^{2,p}$ weak topology to an immersion $\psi:M\to\Ri$.\\
With the arguments of the proof of Theorem~\ref{landel}
in~\cite{della1,langer2} and keeping into account that in our case we
have also the estimates~\eqref{bbbbb}, it is possible to conclude that
actually the convergence is in the $C^{\infty}$ topology and the limit
hypersurface is smooth (see also~\cite{huisk3}, Prop.~3.4).

\begin{teo}
The family of smooth hypersurfaces $\varphi_0:M\to\Ri$, immersed 
in $\Ri$, evolving by the gradient flow for the functional 
$$
{\mathcal F}_m(\varphi)=\int_M1+\vert\nabla^m\nu\vert^2\,d\mu\,,
$$
when $m>\left[\frac{n}{2}\right]$, up to reparametrizations and
translations, is compact in the $C^\infty$ topology of maps. Moreover,
every limit point for $t\to+\infty$ is a $C^{\infty}$ critical
hypersurface of the functional ${\mathcal F}_m$.
\end{teo}

\section{Some Remarks and Open Problems}\label{remmmm}

\subsection{Other Ambient Spaces}

A natural extension would be to consider an ambient
spaces different by $\Ri$ and a codimension $s$ greater than one, 
that is, a general Riemannian manifold $(N,h)$ of dimension $n+s$ 
(notice that Polden's Theorem~\ref{smalltime} about small time
existence of the flow already deals with hypersurfaces in a general
target manifold).
In this context a functional which could be considered is 
$$
{\mathcal F}_m(\varphi)=\int_M1+\vert\nabla^m\omega\vert^2\,d\mu
$$
where $\omega=\nu_1 \wedge \dots \wedge \nu_s$ is a $s$--vector obtained
by a local orthonormal basis of the normal space to the
$n$--dimensional immersed submanifold $\varphi:M\to N^{n+s}$.\\
In~\cite{kuschat1} Kuwert and Sch\"atzle announce a forthcoming paper
with the extension of Polden's results to space curves.

\subsection{Other Functionals}

It would be very interesting to study the flows in the ``critical''
case $2m=n$, where our proof fails since we are no more able to bound
the constants independently of time. 
Notice that the well known Willmore functional 
(see~\cite{kuschat1,simon2,willmore}) 
$$
{\mathcal W}(\varphi)=\int_M \vert\AAA\vert^2\,d\mu  
$$
falls exactly in this case if we add the area term, since $\vert
\AAA\vert^2$ is equal to $\vert\nabla\nu\vert^2$.\\
To the author knowledge, up to now nor there is a proof of regularity of the
flow, neither an example showing the development of a singularity. A
first step in this research was recently done by Kuwert and
Sch\"atzle~\cite{kuschat1}.

When $2m<n$ we do not expect regularity of the flow by the gradient of 
${\mathcal  F}_m$ since, by analogy with the previous discussion about 
the regularity of varifolds, the curvature term should not be sufficient
to give regularity and dumb--bell like separation phenomena should
appear during the flow of certain hypersurfaces. It should also be
noticed that in this and in the critical case, the $n$--dimensional
unit sphere in $\Ri$ collapses in finite time.

Moreover, one can consider also ``non--quadratic'' functionals, 
for instance,
$$
{\mathcal F}_{m,p}(\varphi)=\int_M 1 +\vert
\nabla^m\nu\vert^p\,d\mu\qquad\text{ when $mp>n$}
$$
(following the analogy with the Sobolev spaces), in particular,
$$
{\mathcal F}_{1,p}(\varphi)=\int_M 1 +\vert \AAA \vert^p\,d\mu\qquad\text{ 
  for $p>n$}
$$
which would give rise to a flow of order lower than the one of
${\mathcal F}_m$ when $n>1$.\\
In the same spirit another interesting functional is 
$$
{\mathcal H}_p(\varphi)=\int_M 1 +\vert \HHH \vert^p\,d\mu\qquad\text{ 
  for $p>n$}\,.
$$
In these cases the smoothness of the associated flows is an open
problem.\\

\subsection{Smoothing Terms}

From our analysis, it easily follows that for every positive
constants $\alpha$ and $\beta$ also the gradient flow of the
functional
$$
{\mathcal F}^{\alpha \beta}_m(\varphi)
=\int_M\alpha+\beta\vert\nabla^m\nu\vert^2\,d\mu
$$
exists and it is smooth for every positive time.\\
Moreover, if we consider a general positive geometric functional
$$
{\mathcal G}(\varphi)=\int_M f(\varphi, g, \AAA, \nu, \,\dots\,,
 \nabla^s\AAA, \nabla^l\nu)\,d\mu\,,
$$
such that $f$ is smooth and has a polynomial growth,
choosing an integer $m$ large enough, the gradient flow of the
perturbed functional with $\varepsilon>0$
$$
{\mathcal G}_m^{\varepsilon}(\varphi)=
{\mathcal G}(\varphi)+\varepsilon{\mathcal F}_m(\varphi)
$$
does not develop singularities. This is achieved choosing $m$ so that
the rescaling order of $\vert\nabla^m\nu\vert^2$ is larger than the
rescaling order of $f(\varphi, g, \AAA, \nu, \,\dots\,, \nabla^s\AAA,
\nabla^l\nu)$, in this way the extra terms coming from ${\mathcal G}$
are well controlled  by the leading term in the first variation of 
$\varepsilon{\mathcal F}_m$ and do not affect long time existence.\\
We say that ${\mathcal F}_m$ is a smoothing term for ${\mathcal G}$.

Once we have a sufficiently general family
of smoothing terms we can study what happens varying the parameters,
in particular when the constant in front of them goes to zero.\\
This program, suggested by De~Giorgi's in~\cite{degio5,degio6},
Sec.~5, can be stated as follows:
given a geometric functional ${\mathcal G}$ defined on submanifolds of
the Euclidean space (or a more general ambient space),
\begin{itemize}
\item find a functional ${\mathcal F}$ such that the perturbed
  functionals ${\mathcal G}_\varepsilon={\mathcal
    G}+\varepsilon{\mathcal F}$ give rise to smooth flows;

\item study what happens when $\varepsilon\to0$, in particular, the
  existence of a limit flow and in this case its relation with the
  gradient flow of ${\mathcal G}$ (if it exists, smooth or singular).
\end{itemize}

Our work shows that the functionals ${\mathcal F}_m$ satisfy the
first point for geometric functionals on hypersurfaces in $\Ri$
${\mathcal G}$ with polynomial growth, provided we choose an order $m$
large enough (depending on ${\mathcal G}$).

Concerning the second point,  
a first step would be to consider the possible
limits when $\varepsilon\to0$ of the flows of ${\int_M 1
  +\varepsilon\vert\nabla^m\nu\vert^2\,d\mu}$ 
when $m>[n/2]$ and their relation with the
mean curvature flow. Even the simplest case of the convergence of the
family of flows of curves associated to the functionals
$$
{\mathcal F}^\varepsilon_1(\gamma)=\int_{\SS^1} 1+\varepsilon k^2\,ds
$$
to the mean curvature flow is an open problem.

\subsection{De~Giorgi's Conjecture}

Finally we introduce the original smoothing terms suggested by
De~Giorgi in~\cite{degio5,degio6}. Given a smooth embedded
hypersurface $M\subset\Ri$, we can consider the squared distance
function $\eta^M(x)=[d(x,M)]^2:\Ri\to\R$ which turns out to be smooth
in a neighborhood of  the hypersurface $M$. Then we define the function
$$
A^M(x)=\frac{\vert x\vert^2-\eta^M(x)}{2}
$$
and its derivatives
$$
A^M_{i_1 \dots i_m}(x)=\frac{\partial^m A^M(x)}{\partial x_i \dots
\partial x_m}
$$
whenever they exist, in particular for every $x\in M$.\\
The quantities $A^M_{i_1 \dots i_m}(x)$ for $x\in M$ are related to
the second fundamental form $\AAA(x)$ of $M$ and to its derivatives up
to the order $m-3$, for instance
$$
\vert A^M_{ijk}(x)\vert^2=\sum_{1\leq i, j, k \leq n+1} [A^M_{ijk}(x)]^2=3\vert
\AAA(x) \vert^2\,.
$$
In general there is a bijective relation between the quantities
$A^M_{ijk}(x)$ and the second fundamental form of $M$ at $x$
(see~\cite{ambman1}). 
In the case of immersed manifold, not necessarily embedded, 
the function $A^M(x)$ can be defined using the property that every 
immersion is locally an embedding.\\
The relations of the distance function with the second fundamental
form make it a valuable tool in the study of the evolution by mean
curvature (see~\cite{ambson,soner1}) and more in general of geometric
functionals and flows (see for instance~\cite{ambman1,delzol1,delzol2}).

De~Giorgi suggested that the gradient flow of the functionals 
$$
{\mathcal {DG}}_m(\varphi)=\int_M 1 +\vert A^M_{i_1 \dots
i_{m}}\vert^2\,d\mu
$$
when $m$ is large enough, does not become singular.\\
By analogy with our work we expect that when
$m>\left[\frac{n}{2}\right]+2$ we  obtain regularity.

The first variations of these functionals has been studied by Ambrosio
and the author in~\cite{ambman1}, Sec.~5.3: the leading term of the
first variation of ${\mathcal {DG}}_m$ turns out to be a constant 
multiple of the leading term of $\EEE_{m-2}$ (see Theorem~\ref{fine})
\begin{equation*}
2m(-1)^{m}\overset{\text{$m-2$~{\rm times}}}
{\overbrace{\Delta\Delta\dots\Delta}} \HHH\,,
\end{equation*}
moreover, the functional ${\mathcal {DG}}_m$ has the same rescaling
  properties of ${\mathcal  F}_{m-2}$.\\
The difficult step in repeating our proof stays in
controlling a priori Sobolev and interpolation constants, or more
  precisely in obtaining inequalities of kind
$$
\Vert A^M_{i_1\dots  i_{k}}\Vert_{L^p(\mu)}\leq C \Vert A^M_{i_1\dots
  i_{k+l}}\Vert_{L^q(\mu)}\,,
$$
since the integrals are done on $M$ but the derivatives are taken
along all the directions of the ambient space $\Ri$.

At this moment  the original conjecture of De~Giorgi remains open.

\subsection{Asymptotic Behavior}

An open problem arising from the discussion of the previous section is
  the question of the uniqueness of the limit hypersurfaces. 
It is also unknown to the author if actually it can happen that the
hypersurface goes to the infinity when $t\to+\infty$.\\
To conclude, we mention the problem of classification of the
  limit points of these flows, or equivalently of the critical
  hypersurfaces of ${\mathcal F}_m$. In his
  work~\cite{polden1} Polden completely classifies the limit curves of
  the flow of the functional~\eqref{pden1}, the analogous
  $n$--dimensional result  seems to be a much more difficult task.

\bibliographystyle{amsplain}
\bibliography{sobgeo}

\end{document}